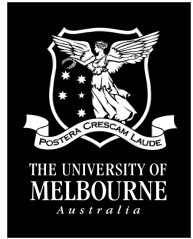

# *Hypercomputation*: computing more than the Turing machine

Toby Ord
Department of Philosophy*
The University of Melbourne
t.ord@pgrad.unimelb.edu.au

**Abstract:**

In this report I provide an introduction to the burgeoning field of hypercomputation – the study of machines that can compute more than Turing machines. I take an extensive survey of many of the key concepts in the field, tying together the disparate ideas and presenting them in a structure which allows comparisons of the many approaches and results. To this I add several new results and draw out some interesting consequences of hypercomputation for several different disciplines.

I begin with a succinct introduction to the classical theory of computation and its place amongst some of the negative results of the 20th Century. I then explain how the Church-Turing Thesis is commonly misunderstood and present new theses which better describe the possible limits on computability. Following this, I introduce ten different hypermachines (including three of my own) and discuss in some depth the manners in which they attain their power and the physical plausibility of each method. I then compare the powers of the different models using a device from recursion theory.

Finally, I examine the implications of hypercomputation to mathematics, physics, computer science and philosophy. Perhaps the most important of these implications is that the negative mathematical results of Gödel, Turing and Chaitin are each dependent upon the nature of physics. This both weakens these results and provides strong links between mathematics and physics. I conclude that hypercomputation is of serious academic interest within many disciplines, opening new possibilities that were previously ignored because of long held misconceptions about the limits of computation.

---

# Acknowledgments

In writing this report, I am indebted to four people. Jack Copeland, for introducing me to hypercomputation and taking my new ideas seriously. Harald Søndergaard, for allowing me the chance to do my Honours thesis on such a fascinating topic and then keeping me focused when I inevitably tried to take on too much. Peter Eckersley, for his ever enthusiastic conversations on the limits of logical thought. And finally, Bernadette Young, for her kindness and support when things looked impossible.

# Preface

It has long been assumed that the Turing machine computes all the functions that are computable in any reasonable sense. It is therefore assumed to be a sufficient model for computability. What then, is one to make of the slow trickle of papers that discuss models which can compute more than the Turing machine? It is perhaps tempting to dismiss such theorising as idle speculation alike to the more fanciful areas of pure mathematics, in which models and abstractions are studied for their own sake with little regard to any real world implications.

This would be a great mistake. The study of more powerful models of computation is of considerable importance, with many far-reaching implications in computer science, mathematics, physics and philosophy. Rarely, if ever, has such an important physical claim about the limits of the universe been so widely accepted from such a weak basis of evidence. The inability of some of the best minds of the century to develop ways in which we could build machines with more computational power than Turing machines is not good evidence that this is impossible. At best, it is reason to consider the problem to be very difficult and likely to require a radically new technique or even new developments in physics.

However, we should not even go this far. A variety of theoretical models for such *hypercomputation* have been presented over the past century. They have, however, often been presented in rather theoretical contexts, which explains the lack of serious work on exploring the physical realisability of these models. Indeed, from what work has been done, there have recently been several serious reports about how quantum mechanics and relativity may be used to harness hypercomputation.

In many ways, the present theory of computation is in a similar position to that of geometry in the 18[th] and 19[th] Centuries. In 1817 one of the most eminent mathematicians, Carl Friedrich Gauss, became convinced that the fifth of Euclid's axioms was independent of the other four and not self evident, so he began to study the consequences of dropping it. However, even Gauss was too fearful of the reactions of the rest of the mathematical community to publish this result. Eight years later, János Bolyai published his independent work on the study of geometry without Euclid's fifth axiom and generated a storm of controversy to last many years.

His work was considered to be in obvious breach of the real-world geometry that the mathematics sought to model and thus an unnecessary and fanciful exercise. However, the new generalised geometry (of which Euclidean geometry is but a special case) gathered a following and led to many new ideas and results. In 1915, Einstein's general theory of relativity suggested that the geometry of our universe is indeed non-Euclidean and was supported by much experimental evidence. The non-

Euclidean geometry of our universe must now be taken into account in both the calculations of theoretical physics and those of a taxi's global positioning system.

Until now, much of the work on hypercomputation has been independently derived, with the occasional new models discussed in relative isolation with only a brief look at their implications. In this report I try to somewhat remedy this situation, presenting a survey of much of the past work on hypercomputation. I explain many of the models that have been presented and analyse their requirements and capabilities before drawing out the considerable implications that they bear for many important areas of mathematics and science.

# Contents

# Chapter 1

# Introduction – Classical Computation

## 1.1   Gödel's Incompleteness Theorem

The desire to capture the complexities of mathematical reasoning and computation in a simple and powerful manner has been a major theme of 20[th] Century mathematics. While developing these ideas, several powerful negative results were reached about the limitations of computation. I refer here to the results of Kurt Gödel, Alan Turing and Gregory Chaitin. Each of these built upon the last to reach more powerful conclusions about the limitations of mathematical knowledge. As well as providing a framework for the rest of this report, these ideas will hopefully be clarified by material in later chapters and will be shown in a new and more positive light.

The 1930's saw two of these great negative results in the foundations of mathematics, the first of which was Gödel's very unexpected proof of 1931. In arguably the most important mathematical paper of the 20[th] Century [27], he proved that the leading formalisation of mathematics, *Principia Mathematica*, was either an incomplete or inconsistent theory of the natural numbers. In other words, that there are propositions in the language of arithmetic that are either true but unprovable within Principia Mathematica or false but provable. Since consistency is required of all serious proof systems, Gödel's result is considered a proof of the incompleteness of Principia Mathematica and is known as *Gödel's Incompleteness Theorem.*

He proved the incompleteness of Principia Mathematica with a combination of two insightful ideas. First, he showed how formulae could be associated with natural numbers, to allow statements in arithmetic to refer to other statements in arithmetic. Then he showed how the predicate 'is provable in Principia Mathematica' could be translated into a statement in arithmetic. This allowed him to construct a statement that says of itself that it is not provable within Principia Mathematica. This statement has the desired property of being either true but unprovable or false but provable.

The existence of such statements for Principia Mathematica would have been considered a major result in itself, but Gödel went further and argued that his method of proof could be applied to *any* formal proof system that operates be finite means. In other words, that all mechanical methods for describing a set of axioms and the rules by which theorems can be inferred can be translated into statements in arithmetic. From these statements, one could then construct new statements that attest to their own unprovability. As a consequence every such system would be incomplete.

This claim that *no* consistent formal system can prove all truths of arithmetic was seen as a profoundly negative result for mathematics. While on the surface it is only a limitation for our formalisation processes (and Gödel himself believed that it was not a limitation for human reasoning [21]), it was taken by many to show that there are statements in arithmetic whose truth is unknowable. Gödel's Incompleteness Theorem is thus usually considered to be a major limitation on the power of reasoning.

## 1.2   Turing Machines and Computation

The second great breakthrough is intimately related to the first – completing it and extending it. Although Gödel's argument was generally accepted by the mathematical community as showing the incompleteness of all formal proof systems, it relied upon the intuitive notion of a generalised formal system. While the proof of the incompleteness of Principia Mathematica was finished and unquestionable, the incompleteness of *all* formal systems was dependent upon a concrete interpretation of a 'mechanical method' or 'finite means'.

The formal analysis of computation by finite means was also motivated by other problems of the day. The *Entscheidungsproblem*[1] and Hilbert's tenth problem[2] both asked for algorithms to solve important mathematical problems. The failure of past efforts to find these algorithms had made some mathematicians suspect that there were, in fact, no algorithms to solve them. Proving the absence of a mechanical method for doing something required a precise formalisation of what mechanical computation is. 1936 saw the independent development of three influential models of computation, aimed at doing just this: the lambda calculus, recursive functions and Turing machines.

All three were soon found to agree on which functions were computable, differing only in how they were to be computed. Alonzo Church's lambda calculus [11] and Steven Kleene's recursive functions [33] were arguably more elegant in this respect, but it was the mechanical action of Turing's machines [52] that most agreed with intuitions about how people go about computing mathematical functions. Unlike the others, it was as much an instruction manual on how such a device might be built as it was a formalism for studying computation. This link between abstract computability and physical computability made Turing machines quickly become the standard model of computation. It also makes Turing machines a natural object for studying even more powerful models of computation, as is done later in this report.

A Turing machine consists of two main parts: a 'tape' for storing the input, output and scratch working, and a set of internal states and transitions detailing how the machine should function.

The tape is made up of an infinite sequence of squares, each of which can have one symbol written on it. It is usually considered as a long strip of paper that is bounded on the left, but stretches infinitely to the right. While this infinite length of the tape may seem quite unrealistic, it is normally justified by considering it to be finite but unbounded – we can imagine that the machine begins with a finite sequence of squares and new blank squares are produced and added to the tape as needed.

---

[1] The problem of deciding whether an arbitrary formula of the predicate calculus is a tautology.

[2] The problem of solving an arbitrary diophantine equation (a polynomial equation with integer variables).

The tape begins with the input inscribed upon it and then, as the computation progresses, this gets gradually altered until it eventually contains just the output. The symbols that can be inscribed on the tape come from a fixed and finite alphabet, $\Sigma$, and each square can hold one symbol. Initially a finite sequence of symbols is placed upon the tape and the rest is left blank. While we could use quite complicated alphabets, we will usually use an alphabet consisting of just two symbols: **0** and **1**. We equate **0** with the blank square and assume the tape to begin with a finite quantity of **1**'s placed at finite distances from the beginning of the tape and the rest is filled with **0**'s.

This tape is linked to the operation of the machine through a 'tape head' which, at any time, is considered to be looking at one particular square of the tape. The instructions of the machine can make this tape head scan the current symbol, change the symbol, move to the left or move to the right.

When the machine is started, it begins in a specified 'initial state'. Associated with this state are a finite list of instructions each of which has four parts (*symbol 1*, *symbol 2*, *direction*, *new state*). The machine treats this instruction as follows: if *symbol 1* is at the tape head then replace it by *symbol 2*, move the tape head in the specified *direction* and then follow the instructions in the *new state*. If an instruction is not applicable (i.e. *symbol 1* is not under the tape head) then it tries the next instruction. At most one instruction can be applicable at a given time and thus there can only be as many instructions for a given state as there are symbols. If no instructions are applicable, the machine halts and the computation is over, with the output being whatever is written on the tape at that time. If a given computation does not halt, it is said to *diverge* and does not produce any output.

While computation here has been defined over functions from finite strings in $\Sigma$ to finite strings in $\Sigma$, it is also possible to consider the computation of other types of functions. For now, I will only examine functions from natural numbers to natural numbers. I will represent the number *n* in unary as a **0** followed by *n* **1**'s (followed by an infinite sequence of **0**'s).

The following is a short program to demonstrate this idea. It simply takes a representation of a natural number and returns 1 if it is even and 0 if it is not. Figure 1 is represents this Turing machine graphically.

state 0:  ( **0**, **0**, right, 1 )

state 1:  ( **1**, **1**, right, 1)
          ( **0**, **0**, left, 2)

state 2:  ( **1**, **0**, left, 3)
          ( **0**, **0**, right, 4)

state 3:  ( **1**, **0**, left, 2)

state 4:  ( **0**, **1**, left, 6 )

state 5:  [no instructions]

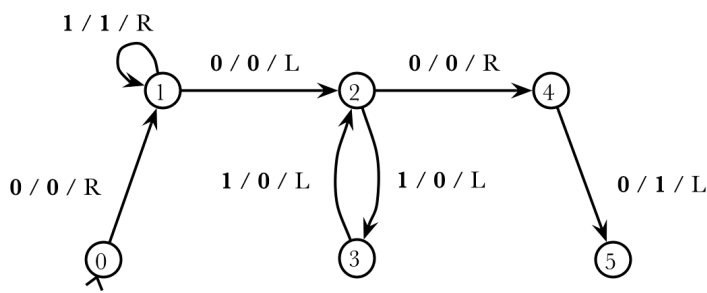

*Figure 1*

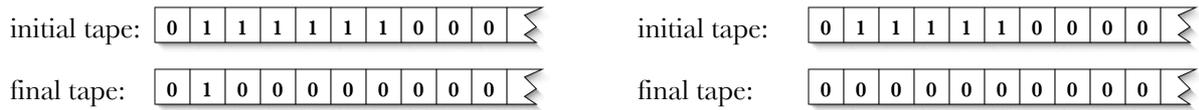

To represent a real number (between 0 and 1) we can use an infinite sequence over {**0**, **1**} to represent the digits of the real in binary[3]. Although a computation that generates such a sequence never finishes, there is a sense in which the real is being computed: if we want to know the value of the real to $n$ bits of accuracy, we can wait until the Turing machine produces the first $n$ bits.

There is a slight catch however, which is that when computing a complicated real such as π, we need some way of knowing which digits are finalised and which are not. To do this, we can give the machine another symbol **#** such that it prints **#** between every two bits of the real and never alters anything to the left of a **#**[4].

For example, this Turing machine (starting with a blank tape) computes the value of ⅓ in binary:

state 0:   ( **0**, **0**, right, 1 )

state 1:   ( **0**, **#**, right, 2 )

state 2:   ( **0**, **1**, right, 3 )

state 3:   ( **0**, **#**, right, 0 )

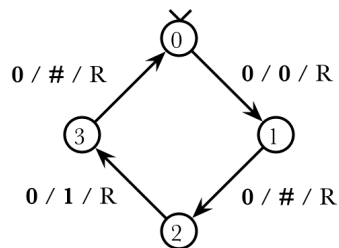

*Figure 2*

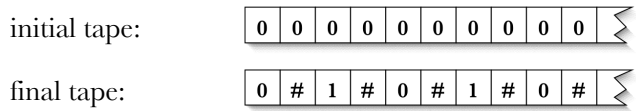

An infinite set of natural numbers can be represented by a real whose $(n+1)^{\text{th}}$ digit is 1 if $n$ is in the set and 0 otherwise. A set of natural numbers is computable if and only if its corresponding real is computable. Using this encoding, the previous Turing machine computes the set of odd naturals.

While they seem quite simple, Turing machines are able to compute many different functions. Indeed, Turing machines were shown to be able to compute anything that the other proposed models could compute and vice versa. Each model was shown to be able to compute the same set of functions and these functions seemed to be equivalent to those that a mathematician could compute without the use of insight.

Many independent models of computation have been produced which also compute the same functions as a Turing machine. Indeed, considerable work has been done on adding extensions to

---

[3] There is no significant problem in representing arbitrary reals, but the arguments are much more clear and concise if attention is restricted to reals between 0 and 1.

[4] Another prominent method for addressing this issue is to give the machine multiple tapes, in particular a read-only input tape, a standard tape for scratch working and a write-only tape for output. These methods are equivalent.

Turing machines – such as allowing them to operate non-deterministically or giving them multiple tapes with their own tape heads – but the enhanced models still computed the same functions. This lent further importance to this set of functions which have come to be known as the *(partial) recursive* functions and even the *computable* functions.

## 1.3   Turing Machines and Semi-computation

It was, however, known since the inception of Turing machines that there were functions that they could not compute. This is clear when you consider that the set of functions from **N** to **N** is uncountable while the set of Turing machines is countable. Both Turing and Church went beyond this however, and showed specific examples of functions which could not be computed [52, 11]. The most famous of these is Turing's *halting function* which takes a natural number representing a Turing machine and its input, returning 1 if the machine halts on its input and 0 if it does not. From this, we also get the existence of a specific uncomputable set $\{n \mid n$ represents a Turing machine / input pair that halts$\}$ and an uncomputable real (the $n^{\text{th}}$ digit of which is 1 if $n$ represents a Turing machine / input pair that halts and 0 otherwise). This set is known as the *halting set* and the real number as τ, for Turing's constant[5].

The halting function can be shown to be uncomputable by Turing machines through a *reductio* argument. Assuming there is a Turing machine for the halting function, $H$, we can construct a modified machine $H_1$ which takes an encoding of a Turing machine and determines whether this machine halts when given its own encoding as input. Now consider another machine $H_2$ which is like $H_1$, but loops if the given machine halts on its own encoding and halts if the given machine loops on its encoding. If $H_2$ is given an encoding of itself as input, it loops if and only if it would halt. This is a contradiction. Since all other claims are verifiable, our assumption of the existence of $H$ must be in error.

While the halting function may not be computable by Turing machines in the terms defined previously, there is still a sense in which it is partially computable. We can construct a Turing machine, $H'$, that takes a representation of a Turing machine, $M$, with some input, $I$, and simulates the behaviour of $M$ on $I$. If $M$ does not halt, then $H'$ does not halt either. If $M$ halts, then $H'$ also halts, but instead of outputting the result of $M$ on $I$ it just outputs the number 1. Thus, $H'$ computes the halting function correctly if the value of the function is 1 and diverges otherwise. It is said that $H'$ *semi-computes* the halting function.

In general, we can say a function, $f : X \to \{0, 1\}$ is *semi-computable* by Turing machines if there is a Turing machine which transforms input, $x$, to $f(x)$ whenever $f(x) = 1$ and either returns 0 or diverges when $f(x) = 0$. Such functions are also said to be *recursively enumerable*.

Thus, the halting function is semi-computable by Turing machines.

---

[5] There are many trivial 1-1 mappings from the set of Turing Machine / input pairs to the set of natural numbers, and for these purposes it does not matter which one is used. The different mappings give different halting functions (and different τ's) but these all have the same important properties.

We also say that a function, *f : X -> {0, 1}* is *co-semi-computable* by Turing machines if there is a Turing machine which transforms input, *x*, to *f(x)* whenever *f(x)* = 0 and either returns 1 or diverges when *f(x)* = 1. Such functions are also said to be *co-recursively enumerable*.

Thus the looping function (does a given Turing machine *not* halt on a given input) is co-semi-computable by Turing machines.

This concept of semi-computability also applies to Turing machines which compute sets or reals. For reals, being semi-computable corresponds to not knowing at what time each digit becomes correct. For example, one could compute τ by letting all the digits start as **0** and simulating all the Turing machine/input combinations in a dovetailed manner, setting the $n^{th}$ bit to **1** when the $n^{th}$ combination halts. In this manner, for each bit that would eventually become **1**, there is a time-step at which this would occur. In other words, if we interpret the tape at any time-step as a real, the machine semi-computes the real *r* if and only if the reals at each time-step converge to *r* from below.

The co-semi-computable reals are those whose complements are semi-computable (where the complement of *r* is defined as 1 - *r*). This is equivalent to starting all the digits as **1**'s and gradually changing some of them to **0**, computing the limit of *r* from above.

As before, we can represent sets of natural numbers with reals so that a set is (co)-semi-computable if and only if its representative real is (co)-semi-computable. Intuitively, this means that for each member of a semi-computable set, there is a time at which it is determined to be in the set and for each natural number that is not a member of a given co-semi-computable set, there is a time at which it is known to not be in the set.

The Turing machine is exactly the type of model that Gödel was looking for to represent his arbitrary formal systems. In a postscriptum to his lectures of 1934 [28], Gödel praises Turing's Machines as allowing a 'precise and unquestionably adequate definition of the general concept of [a] formal system'. With the help of Turing's work on computability, formal systems can be specified as Turing machines that semi-compute a set of formulae, which are considered proven. This can be considered in the terms of classical proof procedures as a recursively enumerable set of axioms with recursively enumerable rules of inference. Gödel's Incompleteness Theorem can therefore be completely specified, stating that no consistent formal system of this type can prove all truths of arithmetic, or, that the set of true formulae of arithmetic is not recursively enumerable.

Turing's proof of the uncomputability of the halting function by his machines also extended Gödel's Incompleteness Theorem. Turing (and Church) had shown an 'absolutely' undecidable function whose values could be proven by no consistent formal system. In contrast, Gödel's previous result only showed that each consistent formal system had its own unprovable result – there was still some hope that each statement of arithmetic was provable in *some* formal system.

Turing's results were thus quite mixed in their portent. On the one hand, Turing gave a theory of computation that provided a deep understanding of algorithmic procedures and, of course, ushered in the current period of digital computing. On the other hand, he provided another alarmingly negative result for the foundations of mathematics, by showing that undecidability in mathematics was even more widespread than had been demonstrated before.

## 1.4 Algorithmic Information Theory

The formalised notion of computation above has also been used to explicate the concept of *information*. Intuitively, if we look at strings of **0**'s and **1**'s, some are more complex than others. For example, **0101010101010101** seems much more ordered than **0001010011101100101**. We would also say that the second string is more complex and that it appears more 'random'. If told that one of these strings was generated by the tosses of a fair coin, we would tend to guess that it was the second.

The field of *Algorithmic Information Theory* (founded by Andrei Kolmogorov [35] and Gregory Chaitin [8, 9]) explains our intuitions by associating the complexity of a string with the size of the smallest program that generates it. Formally, a function $H$ is defined which represents the 'algorithmic information content' or 'algorithmic entropy' of a string, $s$. $H(s)$ is defined as the size of the shortest program (in terms of bits) which produces the string $s$ when given no input[6].

While this definition is relative to the programming language, some very general results are possible. This is because for any two languages ($A$ and $B$) that both compute exactly the recursive functions, there is a fixed length program in $A$ that will act as an interpreter from programs in $B$ into programs in $A$. Thus, across all languages the algorithmic information content of a given string differs by only a constant amount.

With the concept of algorithmic information content in mind, we can say that a string is quite random if the length of its minimal program approaches the length of the string itself. With regards to the two strings above, we would expect the minimal program to generate the second string (in most languages) to be longer than the minimal program to generate the first string. This means that we would expect the first string to be quite compressible in most languages – if we had to send this string to someone using as few bits as possible, we could send the program to generate it. This trick would not work for the second string as its long minimal program (in most languages) makes it effectively incompressible.

While this definition of randomness for bitstrings as length of minimal programs or compressibility seems promising, it suffers in practice since determining the minimal program for an arbitrary string is equivalent to solving the halting problem and is thus uncomputable by Turing machines. We can use a simple counting argument to see that most strings are random, but it is very difficult to show this for a given string.

The concepts of algorithmic information theory extend to infinite strings. For recursive infinite strings, such as the binary expansion of $\pi$, the algorithmic information content is simply the size of the smallest program generating $\pi$. For a non-recursive string such as the binary expansion of $\tau$, there is no finite program, so the algorithmic information content is infinite.

There is, however, a sense in which other non-recursive bitstrings could be more densely filled with information than $\tau$. Consider the case of trying to find out $n$ bits of $\tau$ (i.e., find out which of a group of $n$ computations halt). We could do this recursively if we were told how many of these computations halted. With this information, we could simply run them all in parallel and wait until the specified quantity halt. At this time, we know that the remainder will loop and thus know whether each computation loops or halts, giving us the required *n* bits of τ. To get this information on how many of the programs halt, we only require log(*n*) bits, because these are sufficient to express any number less than or equal to *n*. Thus, the infinite amount of information in τ is distributed very sparsely.

---
[6] Like the halting function, the exact specification of $H$ is relative to the model of computation and the means of encoding programs, but the important properties hold for all (Turing equivalent) models and (recursive) encodings.

Chaitin proposed the constant, Ω, as an example of how an infinite amount of information could be distributed much more densely. Ω is called the 'halting probability' and can be considered as the chance that a program will halt if its bits are chosen by flipping a fair coin[7]. Like τ, the value of Ω is dependent upon the language being used, but its interesting properties are the same regardless [9].

$$\Omega = \sum_{p \text{ halts}} 2^{-|p|}$$

Chaitin defines an infinite bitstring, *s*, to be random if and only if the information content of the initial segment, $s_n$, of length *n* eventually becomes and remains greater than *n*. By this definition, Ω is random while τ is not. In terms of compressibility, we can see that all non-recursive strings are incompressible in that they cannot be derived from a finite program, but only random strings like Ω have the property that of all their prefixes, only a finite amount are compressible.

Continuing the analogy with signal sending to these infinite bitstrings, we see that π can be sent as a finite bitstring by sending a program that computes it. On the other hand, if someone were to have access to the digits of τ and wished to send them to a friend, they would need to send an infinite number of bits. For each *n* bits they send, however, their colleague could construct $2^n$ bits of τ using the method described earlier. τ is thus compressible in a local way, but globally incompressible. If someone had access to Ω, they would need to send *n* bits for every *n* bits their colleague were to receive. Ω is thus incompressible in both ways.

Chaitin has gone even further [9] and translated Ω into an exponential diophantine equation with a parameter, *n*. This equation has the property of having an infinite number of solutions if and only if the $n^{\text{th}}$ bit of Ω is a 1. This moves his result into the domain of arithmetic, along with Gödel's Incompleteness Theorem. It shows that even in arithmetic, there are sets of equations whose solutions have 'no pattern' in the same manner as Ω.

Chaitin has also shown [9] that no consistent formal system of *k* bits in size, can prove whether or not this diophantine equation has infinitely many solutions for more than *k* different values of *n*. The formal system cannot get out more information about the patterns of solutions in this equation that it begins with in its own axioms and rules of inference. Indeed, in predicting this additional information about the patterns of the solutions to this equation, no formal system can do better than chance. In this sense, Chaitin's result can be seen as extending Gödel and Turing's results to not only include incompleteness in the heart of mathematics, but also *randomness*.

---

[7] For this to give a well defined number, we must add a special constraint on the programming language used to generate Ω which is that no program (represented as a bitstring) can be a prefix of another.

In the remainder of this report, I expand upon the ideas of this chapter and present a somewhat different view of these three famous results. This different approach will focus on an element of these results which is often overlooked: the concept of a formal system as a Turing machine. I will explore models of computation that go beyond the power of the Turing machine and examine the interesting and important effects that these have for the questions about bounds of mathematical reasoning and computation.

# Chapter 2

# Hypercomputation

## 2.1 The Church-Turing Thesis

As mentioned earlier, the primary application of the Turing machine model was as a formalised replacement for the intuitive notion of computability. By equating Turing machines and effective procedures, Turing could argue that the absence of a Turing machine to solve the decision problem for the predicate calculus implied the absence of any effective procedure.

It seems quite plausible that all functions that we would naturally regard as computable can be computed by a Turing machine, but it is important to see that this is a substantive claim and should be identified as such. This claim is known as the Church-Turing Thesis [11, 52] and is essential to proofs that certain mathematical functions/sets/reals are uncomputable. While Turing argues convincingly for this thesis, it is not something that can be proven because the concept of an effective procedure is inherently imprecise.

Turing describes an effective procedure as one that could be performed by an idealised, infinitely patient mathematician working with an unlimited supply of paper and pencils – but without insight[8]. It is this that he says is equivalent in power to the Turing machine. He does not deny that other models of computation could compute things that Turing machines cannot. He does not even deny that some physically realisable machine could exceed the power of a Turing machine. He simply equates the power of Turing machines with the rather imprecise notion of a methodical mathematician.

I stress these points because the literature includes many, many, examples of using the Church-Turing Thesis to make these stronger claims. Indeed, it is a rare account of the Church-Turing Thesis that draws the distinction between these possible claims. Jack Copeland's 'The Broad Concept of Computation' [13] includes a valuable account of many of the prominent papers and books that have drastically misstated the Church-Turing Thesis and examines some of the effects of this on the computational theory of mind. In this report, I attempt to rectify some of the confusion that has been caused by this mistake in the fields of mathematics and computer science.

---

[8] Turing's later paper 'Computing Machinery and Intelligence' [54] can be taken as saying that even mathematicians working *with* insight cannot exceed the power of Turing machines.

## 2.2 Other Relevant Theses

Here are a collection of claims about the relative power of Turing machines:

> **A**: All processes performable by idealised mathematicians are simulable by Turing machines
> **B**: All mathematically harnessable processes of the universe are simulable by Turing machines
> **C**: All physically harnessable processes of the universe are simulable by Turing machines
> **D**: All processes of the universe are simulable by Turing machines
> **E**: All formalisable processes are simulable by Turing machines

**A** is the Church-Turing Thesis. **B**, **C**, **D**, **E** are claims that are sometimes confused with it. As mentioned earlier, there are good reasons to believe **A** is true. On the other hand, in the remainder of this report I will be discussing formal processes that are not simulable by Turing machines, thus showing **E** to be false. These models of computation that compute more than Turing machines have come to be known as *hypermachines* and the study of them as *hypercomputation*.

**B**, **C** and **D** are more contentious than the other theses and qualitatively different. Rather than being claims in philosophy or mathematics, they are assertions about the nature of physics and their truth or falsity rests in the underlying structure of our universe.

The distinctions between **B**, **C** and **D** are important ones. By the term *mathematically harnessable* I intuitively mean processes that could be used to compute given mathematical functions. By *physically harnessable* I mean processes that could be used to simulate other given physical processes. Somewhat more formally, I describe the terms as follows:

A process, *P*, is *mathematically harnessable* if and only if the input/output behaviour of *P*:
- is either deterministic or approximately deterministic with the chance/size of error able to be reduced to an arbitrarily small amount
- has a finite mathematical definition which can come to be known through science

A process, *P*, is *physically harnessable* if and only if the input/output behaviour of *P*:
- can be scientifically shown to be usable in the simulation of some other specific process

Like the Church-Turing Thesis, these definitions are still somewhat vague. Unlike the Church-Turing Thesis, however, they refer to physical, epistemological and computational properties of processes that seem important to people seeking to clarify the computational limits of the universe. While these definitions require further refining before they can exactly capture the intended concepts, they seem much more amenable to such elucidation than does the Church-Turing Thesis.

According to these definitions, the physical processes that have been historically used for computation such as mechanical relay switching and current flow through transistors are mathematically and physically harnessable. These are both examples of processes that are *approximately deterministic* in the sense used above – quantum mechanics claims that they are really made up of many tiny random processes and that even on a large scale there is a non-zero chance of them behaving in a non-classical manner. The deterministic behaviour that they approximate is however finitely mathematically defined and we came to know about it through science. Thus, these processes can be used to compute trusted mathematical results and simulate physical systems.

It is logically possible that the mathematically harnessable, physically harnessable and non-harnessable processes differ in the computational power needed to simulate them. For example, consider a universe where arbitrarily fine measurements can be made. If there exists a naturally occurring rod of some strange mineral whose length is a non-recursive real, then the process of measuring this length (and producing this real in, say, binary notation) is not simulable by Turing machines. It may be however, that we cannot know anything else about this length (perhaps we cannot even prove it is non-recursive). In this case, we have a non-recursive process that cannot be harnessed for mathematical or physical purposes.

We could further imagine that, while being unable to know a finite mathematical description of the real (such as those of $\pi$ or $\tau$), our best scientific theories may tell us that all other naturally occurring rods of this mineral must also be of the same length. In this case, the measuring of this rod can be used to simulate physical scenarios involving other such rods. Here we have a process that is physically harnessable, but not mathematically harnessable.

The different possible combinations of theses **B**, **C** and **D** present markedly different types of universe that we may live in:

**B C D**   There is no hypercomputation (harnessable or otherwise) in the universe. Turing machines suffice to simulate all processes of the universe and the universe is thus inherently simulable by our current computers. In this case hypercomputation is an interesting theoretical concept with no direct physical analogue in a similar way to non-Euclidean geometry in a Newtonian universe. Hypercomputation may still be studied as an informative generalisation of computability and its theoretical results may allow unexpected developments in classical computing, but it will be considered at least as abstract an idea as that of complex numbers.

**B C**   The universe is hypercomputational, but no more power can be harnessed than that of a Turing machine. The universe does not have enough harnessable power to build a machine that simulates arbitrary physical scenarios. This would be a very negative result for physics.

**B**   The universe is hypercomputational and it is at least theoretically possible to build a hypermachine to simulate physical scenarios that are not simulable by Turing machines. The hypermachines cannot be understood well enough, however, to be usable to compute specific, predefined, non-recursive mathematical structures. The hypermachines may or may not be capable of simulating all arbitrary physical scenarios (we know the physically harnessable and non-harnessable processes are not simulable by a Turing machine, but it does not follow that they are equally difficult to simulate).

[none]   Hypercomputation exists in the universe and is accessible. As above, the universe may or may not be capable of simulating arbitrary physical scenarios, but it is at least theoretically possible to construct hypermachines to compute some predefined non-recursive functions and to simulate some non-recursive physical scenarios. This would be a very positive universe for mathematicians, with many more mathematical theorems provable than is currently believed possible. The decision problem for the predicate calculus may be solvable and the halting function may be computable.

Unlike the standard Church-Turing Thesis, theses **B**, **C** and **D** present concrete constraints that our universe may or may not obey. They divide the different possibilities for the nature of our universe into a set of interesting and intriguing cases. If we were to find out that the Church-Turing Thesis was true (or that it was false), it is not at all clear what this would say about our universe. It would not say anything about the types of computation that are achievable through physical machines or even what is computable by mathematicians through insight.

While there is still considerable vagueness in the concepts that my proposed theses rely on (such as harnessability and simulability) this is of a different sort to the vagueness of the Church-Turing Thesis. The terms I have used really seem to relate to important physical limits on different types of computational power and appear to be susceptible to further clarification.

The theses discussed above are by no means the only physically important claims regarding the limits of computation. While they stratify the possible laws of physics in terms of the accessibility of hypercomputation, it is also possible to look at the degree of hypercomputation achievable in each case, replacing the references to Turing machines with other more powerful models.

Several people have investigated other methods of analysing physical limits of computation. Kreisel suggests *Thesis M* (for mechanical) which states that 'the behaviour of any discrete physical system evolving according to mechanical laws is recursive' [36]. This claim is similar to thesis **D**, but somewhat weaker in its restriction to *discrete* physical systems. This prevents it from ruling out the existence of some physical machine that uses continuous phenomena to compute non-recursive functions.

This formulation is also lacking in not taking stochastic processes into account. A device that can predict any bit of $\tau$ with 51% accuracy would be a usable hypermachine and contradict theses **B**, **C** and **D** above, but not Thesis *M*. Kreisel accounts for this possibility with his *Thesis P* (for probabilistic) which states that 'any possible behaviour of a discrete physical system (according to present day physical theory) is recursive'[9] [36]. However, this thesis is still restricted to discrete systems and also is conspicuously not about the nature of the universe, but about the nature of our present day physical theory.

Deutsch takes the notions of possible outcomes further, linking probability into his very conception of what it is to compute [22]. He then presents a transformation of the Church-Turing Thesis to a physical principle – *the Church-Turing Principle*: 'Every finitely realisable physical system can be perfectly simulated by a universal model computing machine operating by finite means'. Again, this is similar to **D**, but is not restricted to a mere Turing machine. Deutsch is particularly interested in the case where this computing machine is itself physically realisable, in which case it asserts that a single machine can be built which can simulate any physically realisable system. This claim is of considerable interest for physicists and philosophers, stating that the universe has the power to simulate itself.

---

[9] Some additional problems are caused by Kreisel's definition of a possible behaviour as one with non-zero probability. In an infinite sequence of fair coin throws, all outcomes are of zero probability and his thesis thus ignores them. A better way to express his sentiment towards exceptionally unlikely cases may be to say that all possible behaviours *excepting a set of measure zero* are recursive.

There have been numerous other attempts at adapting the Church-Turing Thesis into a more concrete claim (see [39] for examples). This is an interesting process and one which clarifies the questions that are raised by the study of hypercomputation. Making physical versions of the original claims about effective procedures draws out the important issues at stake: what mathematical functions can one compute in the universe? How computationally difficult is it to simulate the universe? Could there be non-recursive but unharnessable processes? We have a very long way to go in answering these questions, but at least we are beginning to ask the right questions and recognise that they need to be answered.

## 2.3 Recursion Theory

The analysis of hypercomputation began in Turing's 1939 paper, 'Systems of Logic Based on Ordinals' [53]. In this, Turing introduced the O-machine which could compute functions that were beyond the power of Turing machines. It was, however, much less practical in its design and not intended to represent a method that could be followed by human mathematicians. Nevertheless, its importance as an abstract tool for analysing and extending the concept of computation was recognised and it has made a great impact in the field of recursion theory.

Recursion theory (or recursive function theory) is the study of the computation of functions from $\mathbf{N}$ to $\mathbf{N}$. The classical models of Turing, Kleene and Church form the basis of this study by defining the sets of recursive and recursively enumerable functions. This is then extended through relative computation with models such as the O-machine, where additional resources (such as oracles) or new primitive functions can expand the set of 'computable' functions. Many ways to explore the relationships between these classes of functions have been found, producing an extensive body of theory [39, 41].

The main approach to this analysis is in terms of the relative difficulty of the functions. Just as multiplication can be reduced to a procedure that just involves addition, so can other non-recursive functions be reduced to each other. Where the reduction of multiplication to addition lets us say that multiplication is not harder than addition, so the reductions between other functions let us show 'not harder than' relations amongst them. When combined with other proofs which show that one function is definitely harder than another (such as the halting function compared to addition), a complex structure emerges and is the focus of much study.

In contrast to this detailed analysis of the mathematical relationships between the functions, there has been only a slow trickle of research into theoretical machines that can compute these functions. These machines show how very different types of resources such as infinite computation length or fair non-determinism affect computational power.

This report examines many of the hypermachines of the last six decades and attempts to unify some of their features. In chapter 3, I give brief introductions to the different models, explaining their individual approaches. In chapter 4, I look at the resources they use to achieve hypercomputation and provide some assessment of their logical and physical feasibility. In chapter 5, I use some of the results of recursion theory to examine the differing capabilities of these hypermachines. In chapter 6, I explore some of the implications that hypercomputation has for a diverse range of computationally influenced fields, most notably mathematics, computer science and physics.

## 2.4 A Note on Terminology

Because it is generally considered by computer scientists that everything that is computable is computable by a Turing machine, those who study hypercomputation have often been forced into some unfortunate terminology, such as 'Computing the Uncomputable'. To rectify this, I make strong distinctions between terms that are often used synonymously.

In this text, I shall use the terms *computable*, *semi-computable*, *decidable* and *semi-decidable* as relative to a given model of computation and the terms *recursive* and *recursively enumerable* as designating specific sets of mathematical objects.

Thus, for Turing machines, the computable and semi-computable functions coincide with the recursive and recursively enumerable functions, but this need not always be the case.

# Chapter 3

# Hypermachines

## 3.1 O-Machines

The paradigm hypermachine is the O-machine, proposed by Turing in 1939 [53]. It is a Turing machine equipped with an oracle that is capable of answering questions about the membership of a specific set of natural numbers. The machine is also equipped with three special states – the call state, the 1-state and the 0-state – along with a special marker symbol, **μ**. To use its oracle, the machine must first write **μ** on two squares of the tape, then enter the call state. This sends a query to the oracle and the machine ends up in the 1-state if the number of tape squares between the **μ** symbols is an element of the oracle set and ends up in the 0-state otherwise.

It is easy to see that an O-machine with the halting set as its oracle can compute the halting function. While this is a trivial example (simply computing the characteristic function of its oracle), many other functions of interest can be computed using the halting oracle. Indeed, this allows all recursively enumerable functions to be computed.

By definition, a function $f$ is recursively enumerable if and only if there is a Turing machine, $m$, such that for all $n \in \mathbf{N}$, $m$ will output 1 if $f(n) = 1$ and either output 0 or loop otherwise. Therefore, an O-machine with the halting oracle can ask its oracle if $m$ halts on input $n$. If $m$ does not halt, the O-machine can return 0 as $f(n)$ must be 0. If $m$ does halt, then the O-machine can simply simulate $m$ on input $n$ and return its results. As a corollary, all co-recursively enumerable functions are also computable by an O-machine with the halting oracle, by simply performing the previous computation, but returning 0 instead of 1 and vice versa.

Given the power and naturalness of the halting set, it is often used in recursion theory to describe different levels of computability. If a set is computable relative to an O-machine with oracle X, it is said to be *computable in X* or *recursive in X*. The set of Turing machine computable functions can be fitted into this hierarchy by imagining a Turing machine as an O-machine with oracle ∅ (the empty set). No O-machine can solve its own halting function (the proof for Turing machines easily generalises) so there is an obvious sequence of increasingly powerful machines. The generalisation of the halting set to O-machines with oracle X is written as X' where ' is known as the *jump operator*. Thus, the halting function (and all r.e. functions) are said to be recursive in ∅', while the halting function for machines with oracle ∅' is recursive in ∅" and so on.

The O-machines have a very interesting property that any function from **N** to **N** is computable by some O-machine, but no oracle is sufficient to allow O-machines with access to it to compute all these functions. This stems from the fact that an O-machine is a combination of a classical part (the finite set of states and transitions) and an oracle (a possibly infinite set of numbers). There are only countably many classical parts an O-machine could have, while there are uncountably many possible oracles. In one sense, this leads to a trivialisation of generalised computation – if we are allowed to choose the oracle once we know the function, it can be solved by using the oracle as an infinite look-up table. However, the notion of computation by O-machines is interesting and well founded, so long as we are looking at the functions that can be solved with a fixed oracle. This distinction is important for several of the other hypermachines as well.

## 3.2 Turing Machines with Initial Inscriptions

Another way to expand the capabilities of a Turing machine is to allow it to begin with information already on its tape. It is easy to see that a Turing machine with a finite number of symbols initially inscribed on its tape does not exceed the powers of one with a blank tape, because one could simply add some states to the beginning of its program which place the required symbols on its tape and then return the tape head to the start of the tape.

However, if we allow the tape to be initially inscribed with an infinite quantity of symbols, the machine's capabilities are increased. With the digits of $\tau$ inscribed on the odd squares of the tape (leaving the evens for the input, scratch work and output), such a machine could quite easily compute the halting function. The initial inscription is really another form of oracle and has the same computational properties. Compared to the O-machine however, Turing machines with initial inscriptions have the advantages of a natural definition and simplicity of operation, while suffering from the need for an explicitly infinite amount of storage space and initial information.

## 3.3 Coupled Turing Machines

In 'Beyond the Universal Turing Machine' [17], Copeland and Sylvan introduce the *coupled Turing machine*. This is a Turing machine with one or more input channels, providing input to the machine while the computation is in progress. This input could take the form of a symbol from the machine's alphabet being written on the first square of the machine's tape. This square is reserved for the special input and cannot be written over by the tape head. As with the O-machine, the specific sequence of input determines the functions that a coupled Turing machine can perform. For example, if the machine had the digits of $\tau$ being fed into it, one by one, it could compute the halting function and all other recursively enumerable functions. There is the slight complication here of the rate at which these additional bits are added, but it is still clear that if they arrive slowly enough, the machine has time to do its needed calculations before the next bit arrives.

Coupled Turing machines differ from O-machines and Turing machines with initial inscriptions, in that they are finite at any particular time. All they require is a real world process that outputs non-recursive sequences of bits.

## 3.4 Asynchronous Networks of Turing Machines

While networks of communicating Turing machines have been discussed in the classical literature and were found to be equivalent to single Turing machines, this result holds only for synchronous networks. In 'Beyond the Universal Turing Machine' [17], Copeland and Sylvan discuss asynchronous networks, with each Turing machine having a timing function $\Delta_k : \mathbf{N} \to \mathbf{N}$ representing the number of time units between it executing its $n^{th}$ action and executing its $(n+1)^{th}$ action[10].

If two machines operate on the same tape with such timing functions, they can produce non-recursive reals or solve the halting function. This power clearly comes from the timing functions and if these are all recursive functions, the additional power disappears since the networks can once again be simulated with a single Turing machine.

## 3.5 Error Prone Turing Machines

A natural extension of this seemingly unharnessable hypercomputation is the error prone Turing machine which is like a normal Turing machine, but sometimes prints a different symbol to the one intended. For Turing machines with an alphabet of just **0** and **1**, printing the wrong symbol just means printing a **1** where a **0** was intended and vice versa. This erroneous behaviour can be defined by an error function, $e : \mathbf{N} \to \{0, 1\}$ where the machine writes its $n^{th}$ symbol incorrectly if and only if $e(n) = 1$.

It is once again clear that a machine like this could compute the halting function and other non-recursive functions, depending on its error function. In some ways, this is the most plausible of the oracle-type machines discussed so far – indeed, it seems that we may already have machines like this. However, it seems most implausible that their potential hypercomputation could be harnessed to produce useful results.

## 3.6 Probabilistic Turing Machines

Turing machines have also been generalised to include randomness [37]. A *probabilistic Turing machine* can have two applicable transitions from a given state. When in such a state, the machine chooses the transition it will take randomly with equal probability for each transition. When given the same input several times, a probabilistic Turing machine may compute different outputs. For a given input, the machine has a probability distribution of possible outputs. A probabilistic Turing machine is said to compute a function if, for each input, the chance of it giving the correct output is greater than $^1/_2$.

Using this notion of computability, it has been shown that a probabilistic Turing machine computes only the recursive functions [37]. The same is true for computing reals or sets of natural numbers. However, there is still a sense in which probabilistic Turing machines can compute non-recursive output. Consider a machine that simply prints a 0 or a 1 with equal probability and then moves the tape head to the right. Such a machine produces a real. Even though the probability of it producing any given real is zero (and the standard theory would thus say that this machine diverges) if we had such a machine, it would surely be producing some real. Indeed, since there are only countably many recursive reals, this machine would produce a non-recursive real with probability one. The same goes for a machine taking natural numbers as input and randomly adding zero or one to them: it would compute a non-recursive function of its input with probability one.

---

[10] Ties can be resolved in some arbitrary, predefined manner.

Considered in this way, probabilistic Turing machines are quite capable of hypercomputation, but not *harnessable* hypercomputation. In many theoretical uses of probabilistic Turing machines, it is sensible to disregard such non-harnessable computation, but the question of whether or not this degenerate form of hypercomputation is possible is very fundamental to the structure of the universe and should not be ignored.

### 3.7  Infinite State Turing Machines

An infinite state Turing machine is a Turing machine where the set of states is allowed to be infinite. This also implies an infinite amount of transitions, but only a finite quantity leading from a given state.

This gives the Turing machine an infinite program, of which only a finite (but unbounded) amount is used in any given computation. We thus obtain a model that is very similar to the O-machine, but with less distinction between classical and non-classical parts.

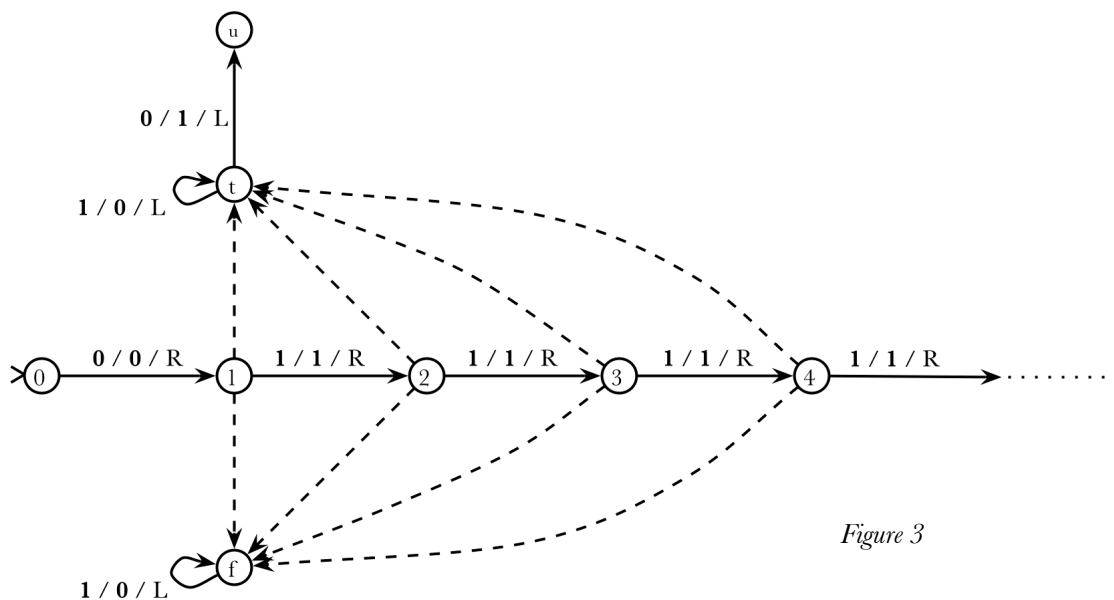

*Figure 3*

Figure 3 shows how an infinite state Turing machine can compute any function from **N** to {0, 1}. It is essentially an infinite lookup table which determines for all numbers given as input, what the result should be. By changing the dotted arcs, any such function can be specified – if $f(n) = 0$, then the dotted arc from state *n* to state *f* should be used (which will replace the number on the tape with the representation of 0), while if $f(n) = 1$, then the dotted arc from state *n* to state *t* should be used (which will replace the number on the tape with the representation of 1). Either way, the transition on the arc should be **1 / 0 / L**.

As can be seen here, the tape is only used for the input and output. If one wanted to define inputs by beginning in different states and outputs by terminating in different states, then the tape could be dispensed with altogether.

## 3.8   Accelerated Turing Machines

In the early 20th Century, Bertrand Russell [42], Ralph Blake [1] and Hermann Weyl [57] independently proposed the idea of a process that performs its first step in one unit of time and each subsequent step in half the time of the step before. Since $1 + {}^1/_2 + {}^1/_4 + {}^1/_8 + \ldots < 2$, such a process could complete an infinity of steps in two time units. The application of this temporal patterning to Turing machines has been discussed briefly by Ian Stewart [48] and in much more depth by Copeland [16] under the name of *accelerated Turing machines*. Since Turing's account of his machines has no mention of how long it takes them to perform an individual step, this acceleration not in conflict with his mathematical conception of a Turing machine.

Consider an accelerated Turing machine, *A*, that was programmed to simulate an arbitrary Turing machine on arbitrary input. If the Turing machine halts on its input, *A* then changes the value of a specified square on its tape (say the first square) from a **0** to a **1**. If the Turing machine does not halt, then *A* leaves the special square as 0. Either way, after 2 time units, the first square on *A*'s tape holds the value of the halting function for this Turing machine and its input.

So far, there has been no difference between an accelerated Turing machine and a standard Turing machine other than the speed at which it operates. In particular, *A* has not solved the halting problem because Turing Machines are defined to output the value on their tape after they halt. In this case, *A* does not halt if its simulated machine does not halt. However, the situation described above suggests a simple change which will allow *A* to solve the halting problem – we consider the machine's output to be whatever is on the first square after two time units.

This model of an accelerated Turing machine only computes functions from **N** to {0,1}, but can be extended to functions from **N** to **N**, by designating the odd squares to be used for the special output, with each of them beginning as **0** and only being changed at most once. In this way, a natural number could be output in a unary representation on the special squares. This could even be extended to allowing real output, where all digits of the real would be written in binary on the special squares after two time units of activity.

## 3.9   Infinite Time Turing Machines

This process of using an infinite computation length can be further extended. In their paper 'Infinite Time Turing Machines' [29], Joel Hamkins and Andy Lewis present a model of a Turing machine that operates for transfinite numbers of steps. We could imagine, for instance, a machine that included an accelerated Turing machine (*M*) as a part. It could initiate *M*'s computation, then after two time units, stop *M*'s movements and reset *M* to its initial state, leaving the tape as it was at the end of the computation. It could then restart *M* with its tape head on the first tape square, running it for another two time units. In such a manner, this machine would perform two infinite sequences of steps in succession. One could even imagine a succession of infinitely many restarts, with *M* performing the whole sequence twice as fast each time, leading to an infinite sequence of infinite sequences of steps.

Perhaps surprisingly, such conceptions of infinite sequences followed by further steps are well founded. The number of steps can be seen as the ordinal numbers:

$0, 1, 2, 3, \ldots \omega, \omega+1, \omega+2, \ldots \omega \cdot 2, \omega \cdot 2+1, \omega \cdot 2+2, \ldots \omega^2, \omega^2+1, \omega^2+2, \ldots \omega^\omega, \ldots$

Here the symbol $\omega$ represents the first transfinite ordinal. It is also a *limit ordinal*, having no immediate predecessor. After an accelerated Turing machine computes for two time units, it has performed $\omega$ steps of computation and if the tape is reused on another computation it has performed $\omega \cdot 2$ steps. The infinite sequence of infinite sequences of steps is denoted by $\omega^2$.

The infinite time Turing machine is a natural extension of the Turing machine to transfinite ordinal times. To determine the configuration of the machine at any successor ordinal time, the new configuration is defined from the old one according to the standard Turing machine rules. At a limit ordinal time, however, the machine's configuration is defined based on all the preceding configurations. The machine goes into a special *limit-state* and each tape square takes a value as follows:

$$\text{square } n \text{ at time } \alpha = \begin{cases} \mathbf{0}, & \text{if the square has settled down to } \mathbf{0} \\ \mathbf{1}, & \text{if the square has settled down to } \mathbf{1} \\ \mathbf{1}, & \text{if the square alternates between } \mathbf{0} \text{ and } \mathbf{1} \text{ unboundedly often} \end{cases}$$

The tape head is placed back on the first square and the machine then continues its computation from this limit-state as it would from any other. As usual, if there is no appropriate step to execute at some point, the machine halts. It can thus perform a finite amount of steps and halt, or an infinite amount of steps and halt, or keep operating through all of the ordinal times and never halt.

Such a machine could compute any recursively enumerable function in $\omega$ steps, by setting the first square on its tape to **0**, then evaluating the function, setting the first square to **1** if *f(n)* = 1. If *f(n)* = 1, then after $\omega$ steps, the first square will hold the value of **1** and if *f(n)* = 0, then after $\omega$ steps, the first square will hold the value **0**. A similar method also computes any of the recursively enumerable reals.

Since infinite time Turing machines can use the entirety of their tapes during their execution, it is natural to define them to accept infinite input (inscribed on, say, the odd squares) and produce infinite output. This allows them a much greater scope in the functions they can compute, but I shall restrict my study here to those infinite time Turing machines that take only finite input, to allow more direct comparison with the other models discussed.

## 3.10 Fair Non-Deterministic Turing Machines

Turing machines can also be generalised to behave non-deterministically [45, 49]. Where a (deterministic) Turing machine always has at most one applicable action in any circumstance, a non-deterministic Turing machine can have many. In such cases, the execution can be thought to branch, trying both possibilities in parallel. Each branch of the computation progresses exactly as it would in a standard Turing machine, except that the output is restricted to being 0 or 1[11]. If at least one branch of the non-deterministic computation halts returning 1, then the computation is considered to return 1. If there is no branch that returns 1, and at least one non-halting branch, then the computation diverges. Otherwise (if all branches halt returning 0) then the machine returns 0.

In this manner, a non-deterministic Turing machine can be considered as using parallel processes (or lucky guesses) to quickly compute its solution. While this seems to lead to significant speedups in computation time, it does not lead to any additional computational power. A non-deterministic Turing machine cannot compute any functions that cannot be computed by a deterministic Turing machine. A deterministic machine can simulate a non-deterministic one by computing each of the branches of computation in a breadth first interleaved manner. However, by restricting this non-determinism to *fair* computations (as done by Edith Spaan, Leen Torenvliet and Peter van Emde Boas [45]), the situation is quite different.

A computation of a non-deterministic Turing machine is said to be *unfair* iff it holds for an infinite suffix of this computation that, if it is in a state $s$ infinitely often, one of the transitions from $s$ is never chosen. A non-deterministic Turing machine is called *fair* if it produces no unfair computations[12].

Spaan, Torenvliet and van Emde Boas show an interesting manner in which a fair non-deterministic Turing machine, $F$, could solve the halting function. $F$ would first go to the end of its input and non-deterministically write an arbitrary natural number, $n$, in unary. This can be done by beginning with 0 and non-deterministically choosing between accepting this number and moving on to the next stage or incrementing the number and repeating the non-deterministic choice. Once $F$ has generated its value for $n$ it can then run a *bounded-time* halting function algorithm, to see if the machine/input combination halts in $n$ steps. If it does, $F$ returns 1. If not, it returns 0. It is clear that if there is a time at which the machine/input combination halts, $F$ will have a finite computation branch that returns the correct answer (1) and $F$ will therefore halt.

What if the given combination does not halt? In this case, the only way $F$ will fail to halt is if it has an infinite computation. The only place one could occur is when the arbitrary number is being generated, but the only way this could not halt is by choosing to increment the number an infinite amount of times, avoiding halting each time. This would be an unfair computation and is thus not possible. Therefore, since all the other branches return 0, the machine will also return 0 which is the correct answer. In this strange way, fair non-deterministic Turing machines can compute the halting function.

---

[11] This is not essential and, like accelerated Turing machines, non-deterministic Turing machines can be thought of as returning arbitrary natural numbers, but the presentation of them as returning only a single bit is by far the norm and is considerably simpler. To make this restriction, the output can be considered as whatever is on the first tape square when the computation halts.

[12] This is just one of several types of 'fairness' that are discussed in the study of non-determinism.

# Chapter 4

# Resources

## 4.1 Assessing the Resources Used to Achieve Hypercomputation

The models discussed in the previous chapter use a wide variety of resources to increase their power over that of the Turing machine. Many people familiar with the Turing machine paradigm and the traditional argument that the power of Turing machines is closed under the addition of new resources (such as extra tapes and non-determinism) will find these new types of resources quite implausible. Indeed, in discussing these resources I show how some of them are logically implausible – having somewhat paradoxical properties – and how some are physically implausible – requiring things that seem impossible according to current physics.

While physical implausibility is important with regards to building one of these devices, it is important to note that our current theories about the nature of physics are far from conclusive. We must be careful not to think of arbitrary measurement and superluminal travel as impossible, but rather as *impossible relative to our best theories*. Just as Newtonian physics was overturned and found to be overly restrictive in some areas (forbidding randomness) and overly lenient in others (allowing arbitrarily fast travel) our current theories may also be quite removed from the truth about physics. We are therefore unable to make claims about certain processes being physically impossible or, indeed, physically possible. Of course it is still very instructive to see which resources are possible or impossible with regards to our current theories.

It must also be noted that this is not a complete analysis of the implications that our current physical theories bear for hypercomputation. Since the resources discussed involve some of the more esoteric areas of modern physics, it would be a major undertaking to give a complete account of the physical possibilities for hypercomputation. Instead, this section should be considered as an informal look at the implications of physics on the hypermachines of chapter 3. I describe many physical phenomena that are important to the realisation of a hypermachine, but do not reach any hard conclusions regarding whether or not such machines can be built.

Even if a certain resource *is* physically impossible, it can still be very useful in understanding the theory of hypercomputation. We may, for example, see that it gives equivalent power to another resource that is physically possible, allowing us to use the impossible resource to help us understand the powers of the possible one. Even if no hypercomputation is physically possible at all, a study of the resources which lead to it is still important for understanding the generalised notion of an algorithm.

## 4.2  Infinite Memory

Accelerated Turing machines and infinite time Turing machines require an infinite amount of memory. Unlike standard Turing machines that only use a finite, yet unbounded, amount of memory, these models can require the whole tape to store their working. Therefore, building such machines requires a means of storing (and retrieving) an infinite bitstring. Unlike the other resources discussed in this chapter, infinite memory is not sufficient for hypercomputation, but it is necessary for some models.

The infinite time Turing machine explicitly requires the storage of an infinite bitstring. After a limit ordinal number of time-steps, it's state is reset to a special limit-state and its head is reset to the first square, but the tape is now filled with an infinite bitstring. This bitstring is then used in further computation. Infinite time Turing machines can be infinitely sensitive to the contents of their tape, so finite approximations can not be made[13].

Accelerated Turing machines seem slightly less dependent on their infinite memory, since the only time they use an infinite amount of memory is when they do not halt. In such cases, there is a finite time-step at which their output is correct and will not change. If at any time after this they were to be denied additional memory and their computation was simply halted by some external means, the result of the computation would be unaffected. There is thus a sense in which there is no difference for these machines between an unbounded amount of memory and an infinite amount. This is analogous to the problem of some standard Turing machines requiring an infinite amount of memory when they do not halt. It is thus somewhat unclear as to whether an accelerated Turing machine really does require the ability to store an infinite bitstring. This may be important for the physical realisation of such a machine because, as discussed below, infinite memory may not be physically possible.

The storage of an infinite bitstring presents no logical problems. Indeed, it is arguably the simplest of the resources discussed in this chapter. We can simply use a real number or an infinite binary sequence as the state for the machine. Physical storage, on the other hand, is much more problematic.

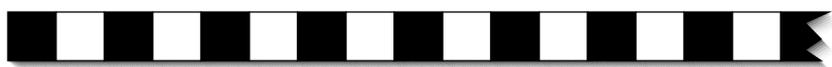

*Figure 4*

An obvious method for storing an infinite bitstring is to take literally the notion of the Turing machine's infinite tape. We could imagine storing an infinite bitstring as an infinite series of white and black cubes floating in space (see figure 4). The obvious disadvantage of such a system is the need for both an infinite volume of space and an infinite amount of matter. Other methods of marking regions of space as **0** or **1** are possible, such as with energy or even different levels of space-time curvature, but they are all somewhat implausible with our current understanding of physics.

---

[13] Consider an infinite time Turing machine that, after its first ω steps, runs through the tape looking for a single **1** and returning 1 if it finds one and 0 if it does not. This machine is 'infinitely sensitive' to the contents of its tape after the first ω steps, because *any* deviation from all **0**'s, no matter how small, will greatly change its output. Thus some infinite time Turing machines have zero error-tolerance.

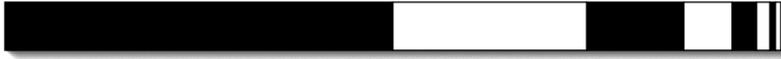
*Figure 5*

Alternatively, we could remove the need for infinite space and resources by marking a bounded area of space with marks of diminishing size (see figure 5). This avoids the problems discussed above, but has the disadvantage of infinite precision being required to store the appropriate bits. According to quantum mechanics, we cannot measure any distance more accurately than the Planck length ($10^{-35}$ meters), so methods like this which require arbitrarily small distances seem impossible.

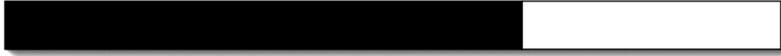
*Figure 6*

A third approach is to use continuous values of other physical quantities (see figure 6). Because the digits of a real can be considered as describing an infinite bitstring, being able to create and measure arbitrary real values would be sufficient for an infinite memory[14]. There are many quantities that we could consider using to store the infinite bitstring, such as the velocity of a particle, the charge on a capacitor or the wavelength of a photon. Unfortunately, it looks as if all such physical variables have limits on how accurately they can be measured [23]. If this is true, then they can only be used to store a finite amount of information.

It is therefore important to note that infinite memory is quite implausible according to our current physical theories. While other exotic means for storing infinite bitstrings have been considered (such as the infinite quantum superpositions suggested in [7]), it is not at all clear that infinite memory is a physical possibility.

## 4.3 Non-Recursive Information Sources

O-machines, Turing machines with initial inscriptions, coupled Turing machines, asynchronous networks of Turing machines, error prone Turing machines and probabilistic Turing machines all have a classical part that is supplemented by some form of non-recursive information source. This constitutes a resource that is similar but importantly different to an infinite amount of memory.

Unlike the infinite memory described previously, non-recursive information sources do not need to be changeable. The machines mentioned in the paragraph above begin with a non-recursive specification, but then only need a finite amount of information to describe their current state. Thus, while these machine need to be able to access an infinite amount of information, they never need to store this much information. However, unlike the contents of an infinite memory which are slowly and routinely built up through simple steps, the initial specifications of these machines consist of non-recursive information from some external source. This is a very important distinction as it raises the questions of where this information would come from and why we should expect it to take a certain form (such as the digits of $\tau$).

---

[14] In contrast, the method of using marks of diminishing size only requires access to *rational* valued intervals.

The machines discussed in this section all achieve hypercomputation through conceptually different means. However, when discussing them from the point of view of physical plausibility, these models are essentially the same. While the different conceptual resources hint at different possibilities for their implementation, the existence of any one of their information sources would allow us to simulate the rest. For instance, if we had a non-recursive initial inscription, we could use it as an oracle for an O-machine or as a timing function for an asynchronous network of Turing machines.

How then might we gain access to a non-recursive information source to create one of the above machines? One way would be to find and use a real valued quantity. If we could find a rod of a non-recursive length, a particle with a non-recursive charge or even a physical constant (such as $c$ or $G$) with a non-recursive value, then the process of measuring the value to greater and greater accuracy would provide another type of non-recursive information source. Since there are uncountably many reals, but only countably many recursive reals, it seems quite plausible that many quantities may take on non-recursive values. Unfortunately, this method suffers from the major drawback (mentioned in section 4.2) of requiring arbitrary precision measurement – this appears to be impossible according to current physical theories.

A better approach may be through the observation of some non-recursive physical process. Consider, for instance, a light that flashes on and off intermittently in a non-recursive, yet deterministic, pattern. Such a light could be used to provide the non-recursive information for any of these equivalent machines and is not in disagreement with current physical theories.

A similar process is actually found in elementary quantum mechanics. Certain particles are unstable and decay after a length of time. Such particles are said to have a half life – the length of time after which the particle has a 50% chance of decaying. Thus, one could design a machine to create a series of such particles and observe each for one half life. If the particle decays within this time, a light could be flashed on, and otherwise left off. According to quantum mechanics, this pattern of flashes would be completely random. Therefore, like the output of the probabilistic Turing machine discussed in section 3.6, this sequence of flashes would be non-recursive with probability one. We could then use the bitstring produced as an oracle for an O-machine and compute non-recursive functions over the natural numbers.

We therefore have a case of quantum mechanics suggesting very strongly that non-recursive processes exist and that the universe is not simulable by Turing machines. However, the plausibility of this non-recursive process comes at the expense of its usefulness. If the pattern of decays is really *completely* random in the sense of probability theory (as distinct from that of algorithmic information theory), then it seems to be a perfect example of a process that is not mathematically harnessable. It may perhaps be physically harnessable (e.g. allowing us to predict the pattern of some other process that has been entangled with it), but it seems to be completely useless for computing definite mathematical results that could not be computed with a Turing machine.

Perhaps some of the other proposed non-recursive information sources may be mathematically harnessable. If the flashes of some light gave us the bits of τ, then we could use this to compute any recursively enumerable function. Is it likely that the bits of τ or some equivalent real would occur in some physical process? On the surface, it would seem not, but there are a couple of ways in which it may be possible. For instance, it may be that there are deep links between physics and computation which could give rise to physical occurrences of such a computationally oriented constant. Alternatively, we may be able to engineer such an information stream – $\tau$ has a very simple specification and we may be able to constrain an otherwise unspecified stream of bits to follow it.

While a natural process giving us the bits of $\tau$ does seem rather implausible, there are some equivalent processes that would be more believable. For example, $\tau$ is known to be reducible to a function known as the *shift function*[15]. This is a function over the natural numbers that is equal to the maximum number of transitions that a Turing machine with *n* states can make before halting. A special property of this function is that no function which is an upper bound for it is recursive. If we had any function that was an upper bound to the shift function, it would be an upper bound to the length of time a program of size *n* would run and would thus allow us compute the halting function. Therefore, if there is any naturally occurring function that is an upper bound to the shift function, we can use it to compute $\tau$. Such a method seems quite promising as it does not require the infinite accuracy of $\tau$: the physical function does not need to have any deep links with theoretical computation, it just needs to grow very, very quickly.

Even if we had some hidden function that could be harnessed like this to produce $\tau$, how could we be sure enough of this to use the results it gives us? We could partially check it (by seeing whether it disagrees with some of the more simple cases that we can calculate ourselves), but we could not completely check that this was really $\tau$, because this would require the same amount of computational resources as are required to generate $\tau$ in the first place. How then could we have any confidence in the accuracy of the input stream?

The only real way is if our best physical theory predicts that this process will produce the bits of $\tau$. In this case, we would have just as much reason to believe that the input stream really was $\tau$ as we have reason to believe that a physical version of a Turing machine really will operate like its mathematical ideal. Indeed, this should give just as much confidence to the accuracy of the input stream as mechanically checking it ever could, because our confidence in any checking process can not be any higher than our confidence in the physical theory that underlies that process. This will remain an additional condition for all of the potential hypercomputational resources: in addition to them allowing hypercomputation, to really harness it, we must be able to scientifically justify that using the resource will give the intended results.

Unharnessable non-recursive hidden functions seem very plausible and even exist in elementary quantum mechanics. However, harnessable non-recursive hidden functions are more difficult to come by. There is little hope of naturally finding a sequence like $\tau$, but more hope of finding an equivalent non-recursive function like the shift function. It is also somewhat plausible that we could engineer a specific non-recursive sequence such as $\tau$.

---

[15] The shift function is intimately related to the more famous (but slightly less intuitive) *busy beaver* function.

## 4.4 Infinite Specification

The models discussed in the previous section have a classical part (which is essentially a Turing machine) and a separate non-classical part that can be provided through an external information source. In contrast, the infinite state Turing machines cannot be divided in such a way. They are therefore naturally thought of as requiring an explicit, infinite specification.

This is highly problematic. The task of providing an arbitrary infinite specification seems far beyond anything we could hope to do. While we can understand certain infinite bitstrings (such as the binary expansion of $\pi$), there are only a countable amount that can be finitely expressed within any given mathematical formalism. To gain the true power of infinite state Turing machines, we would therefore require the ability to provide arbitrary infinite specifications. Such a task is clearly beyond our unaided power. Indeed, since infinite time Turing machines can be constructed to compute each function from $\mathbf{N}$ to $\{0,1\}$ using an infinite lookup table (see section 3.7), specifying an infinite time Turing machine to compute an arbitrary function does not seem any easier than computing the function in the first place.

Infinite state Turing machines are therefore best regarded as theoretical models which show how allowing a Turing machine to use an infinity of states creates a degenerate kind of computer – one that is so powerful as to trivially compute all functions from $\mathbf{N}$ to $\{0,1\}$. This is still an interesting result and may lead to interesting comparisons with other infinitely specified models for computing functions from $\mathbf{R}$ to $\{0,1\}$ (which would *not* be able to trivially compute each of these functions). However, the infinite state Turing machine is a limiting case of discrete computation and it's method for achieving hypercomputation does not seem to be useful for constructing a physical hypermachine.

## 4.5 Infinite Computation

Performing an infinite number of steps in the computation is a very different method of hypercomputation. The accelerated Turing machine already shows a manner in which this could be possible: performing the steps faster and faster to get through an infinite amount in a finite time.

However, to achieve infinite computation through acceleration like this, we rapidly run into conflict with modern physics. For the tape head to get faster and faster, it would need to get faster than the speed of light. This is not quite a breach of the theory of special relativity as objects are allowed to exceed the speed of light so long as they are never slower than it. This is however on the very edge of physical plausibility.

A better approach could be to decrease the distance that the head needs to travel at each step, thus avoiding the need for exceeding the speed of light. Such an approach is adopted by Brian Davies in 'Building Infinite Machines' [19], who also attempts to solve the problem of excessive heat production. This approach does have the drawback, however, of requiring infinite spatial precision which has conflicts with quantum mechanics as I have discussed above.

Mark Hogarth [30, 31] has argued that general relativity allows infinite computation in certain exotic space-times (known as Malament-Hogarth space-times). In such space-times, it may be possible to set up a computer and an observer in such a way that the observer witnesses the computer undertake an infinite number of steps within a finite time. Further examination of such systems has been undertaken by Gábor Etesi and István Németi [25] who pay particular attention to the Malament-Hogarth space-time predicted to occur around a rotating, electrically charged black hole – a situation that they see as quite physically plausible. An in depth examination of Hogarth's particular model of hypercomputation (which differs from those presented here) can be found in Lucian Wischik's 'A Formalisation of Non-Finite Computation' [58].

In their recent paper, 'Coins, Quantum Measurements, and Turing's Barrier' [7], Cristian Calude and Boris Pavlov present a model for using quantum mechanics to achieve infinite computation. Their method is quite different to the accelerated Turing machines discussed here, in that it does not perform the infinite quantity of steps in series, but rather in parallel. This gives them a distinct and interesting way of achieving infinite computation, showing how quantum mechanical properties may help hypercomputation, rather than limit it.

Leaving aside the physical difficulties with such a system, there are also some logical problems. In particular, there is the paradox of Thomson's lamp [51]. If we ran an accelerated Turing machine with a program which simply replaced the first square of the tape with a **1** and then a **0** and repeated this procedure over and over, what would be on this square after two time units? There does not seem to be a good answer. The specification of an accelerated Turing machine made sure that one could only change the designated square(s) once to avoid exactly this type of problem. However, this situation still arises on the standard parts of the tape. While we do not need to look at these squares for the computation to work, it seems quite problematic if they may have no defined value after two time units.

The infinite time Turing machine avoids this by carefully defining the value of a tape square after an infinite number of time-steps to be **0** if the symbol on the square is eventually **0** (i.e. if it becomes **0** at some stage and then doesn't change again) and **1** if it is eventually **1** or if it alternates between **0** and **1** unboundedly often. This overcomes the paradox and gives the machine a well defined behaviour, but makes it more challenging to see how the process could be implemented.

It has also been argued (by Karl Svozil [50]) that accelerated Turing machines are paradoxical for being able to compute the halting function. Svozil argues that the proof that no Turing machine can compute the halting function (because this would allow a machine with contradictory properties to be built) implies that no accelerated Turing machine could compute it either. One way of avoiding this problem would be to say that the original proof only restricts machines from computing their own halting functions and we have no reason to believe that an accelerated Turing machine could do this. However, an accelerated Turing machine is exactly the same as a standard Turing machine except for its method of outputting the result. In particular, their halting behaviour is identical, thus accelerated Turing machines *can* actually compute their own halting functions[16].

---

[16] Note that this is not the case for infinite time Turing machines and this counter argument therefore succeeds in defending infinite time Turing machines from this attack.

The way out of this paradox is that in the original proof, a new machine was created that solves the halting function for a given machine and its input, halting if and only if that machine does not. Such an accelerated Turing machine cannot be constructed. If an accelerated Turing machine is set to compute the halting function and its input specifies a non-halting computation, this machine cannot determine this *and then halt*. Determining that this machine does not halt involves the accelerated Turing machine looping. Thus the paradoxical machine cannot be created and the argument for the contradictory nature of accelerated Turing machines fails[17].

Unlike the accelerated Turing machine, the infinite time Turing machine can perform an infinite amount of computation and then go on to compute more. For example, it could compute a recursively enumerable function by first computing the value of $\tau$ (in an infinite number of steps) and then using this as an oracle to compute the function (in a finite number of additional steps). Indeed, the infinite time Turing machine is in a well defined state at each ordinal time. This presents no logical problem, but may well make it physically impossible.

Restricting the maximum allowable time-steps for an infinite time Turing machine gives a hierarchy of increasingly powerful machines. We can consider a standard Turing machine as an infinite time Turing machine restricted to less than $\omega$ time-steps. We can then define two types of infinite time Turing machine for each limit ordinal $\lambda$, a '$\lambda$-time Turing machine' and a '$<\lambda$-time Turing machine'. Two particularly important results regarding this are Hamkins and Lewis' result that no infinite time Turing machine computing for an uncountable ordinal time can halt [29], and my result in chapter 5 that $<\omega^2$ time-steps suffice for computing all sentences of arithmetic.

We can further separate the powers of infinite time Turing machines by the number of 'designated squares' that they possess. This is just one for an accelerated Turing machine and every square for an infinite time Turing machine, but this need not be the case. We could define Turing machines that compute for more than $\omega$ time-steps, but can only save a finite amount of bits as they pass each limit ordinal time (with the other information being lost).

We could also restrict the amount of times that the symbol on a designated square can be changed. This number is one for accelerated Turing machines and $\omega$ for infinite time Turing machines, but again this need not be the case.

These changes allow for a continuum of possible hypermachines between the power of the accelerated Turing machines and the infinite time Turing machines. This is clearly of interest to the logical study of infinite computation and could be of physical relevance as well, if the physical laws limit computation at a point between the two extreme models.

Despite potential problems with implementation, infinite computation is a very interesting theoretical resource because it does not gain power from having additional information injected into the model, but rather allows the computation to explore every last implication of the initial information.

---

[17]Copeland [16] gives a similar argument, comparing the methods of outputting the results of each of accelerated and standard Turing machines. He distinguishes between the 'internal' computation of the standard Turing machine which allows it to use its results for further computation and the 'external' computation of the accelerated Turing machine which does not.

## 4.6 Fair Non-Determinism

To create a fair non-deterministic Turing machine, we would either need access to some physical process that directly implements fair non-determinism or a process that we could use to enforce fairness on a simulated non-deterministic system. Neither option is denied by our current physical theories, but there is only limited support. The most promising avenue seems to be through quantum mechanics. The standard model of quantum computing (see Deutsch [22]) does not allow the computation of any non-recursive functions, but does use a form of non-determinism to speed up the computation. There may be a way to enforce fairness on such a process, but I am not aware of any method for doing so. The question of whether fair non-deterministic processes can exist according to our current physical theories is very much open.

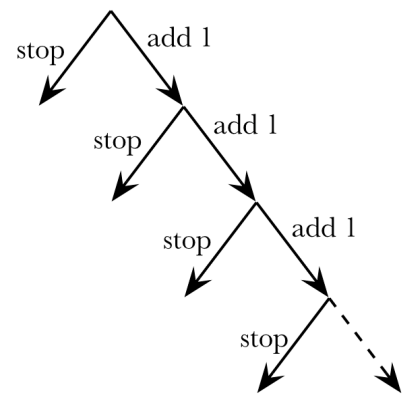

*Figure 7*

On the other hand, there are reasons to be suspicious of the very concept of fair non-determinism. These stem from the fact that fair non-determinism is a form of unbounded non-determinism. This means that a fair non-deterministic Turing machine, such as *F* (described in section 3.10) will always halt, but can do so in an infinite number of different ways (generating any $n \in \mathbf{N}$ in the process). Gordon Plotkin [40] points out a tension between this and König's lemma (which states that if all branches of a finitely branching tree are finite, then so is the tree itself). Plotkin argues that therefore the computation tree of *F* must have an infinite branch and hence the possibility of non-termination (see figure 7 for the computation tree involved in *F*'s generation of *n*). However, as pointed out by William Clinger [12], it is possible that the computation tree does indeed contain an infinite branch which has each of its nodes present in one of the halting computations, but there does not need to be any single computation that proceeds *all* the way down the infinite branch. This counters Plotkin's argument, but shows that fair non-determinism has some very unintuitive consequences.

Another example of the unintuitive nature of fair non-determinism can be seen in the behaviour of *F*. Recall that *F* computes the halting function by generating an arbitrarily large integer, *n*, and testing to see if the given machine halts on its input within *n* time-steps. If the simulated machine halts, then we would imagine that the time taken by *F* is related to the smallest value of *n* that was required. What if the simulated machine does not halt? In this case we know that *F* will return 0, but it is entirely unclear as to how this will happen. This problem becomes particularly noticeable when we consider the length of time *F* will take to return its result. If we consider the non-determinism as computing each branch in parallel, then we would assume that the machine would halt when all the branches have halted, however, even though each branch halts, there is no time at which they have all halted. Although this machine has no infinite branch, its infinite sequence of increasingly long finite branches seem to keep it from halting just as effectively.

The status of fair non-determinism is therefore somewhat unclear. There do not seem to be any strong reasons to believe that it is physically impossible, but its logical consequences are quite unintuitive and need to be better understood before we can be sure that it is free of internal contradiction.

## Chapter 5

# Capabilities

## 5.1 The Arithmetical Hierarchy

Following recursion theory, we can develop a means of assessing the relative powers of different hypermachines by examining the sentences of arithmetic. These are well formed formulae composed of natural numbers, variables ranging over natural numbers, $+$, $\cdot$, $\neg$, $\wedge$, $\vee$, $\forall$, $\exists$, $=$. These formulae include many elementary statements such as:

3 does not divide 8:     $\neg \exists x \, (8 = 3 \cdot x)$
Addition is commutative:  $\forall x \, \forall y \, (x + y = y + x)$

It is also possible to express more complex statements about the natural numbers whose truth is much less clear. Indeed, it is this language for which Gödel proved his incompleteness theorem – no consistent formal system can prove all the true formulae of arithmetic. In his proof [27], Gödel showed how predicates could be constructed in this language through leaving free variables in the formula. For example, we can have:

$x \mid y$:     $\exists z \, (x = y \cdot z)$
*Prime*($x$):  $\neg (x = 1) \wedge \neg \exists z \, ( \neg (z = 1) \wedge \neg (z = x) \wedge x \mid z )$

It is this form of 'programming' in the language of arithmetic that Gödel used to construct the predicate 'is provable in Principia Mathematica'. From this, he created his famous statement that claims it is unprovable within Principia Mathematica. Gödel's generalisation of this argument to find unprovable statements for all formal systems shows that this language of first order arithmetic is very expressive.

The full power of its expressiveness was explored much more fully by Kleene [34] who stratified the predicates expressible in first order arithmetic into what is called the *arithmetical hierarchy*. This consists of the classes $\Sigma_n$, $\Pi_n$, $\Delta_n$, where $n \in \mathbf{N}$. These are defined as follows[18]:

---

[18] Note that I only look at formulae with one free variable (and thus predicates of one variable). This can be done without loss of generality because of the existence of recursive mappings from $\mathbf{N}^m$ to $\mathbf{N}$.

$\Sigma_1$: $\{x : \exists y\, (\, R(x, y)\, )\}$      where $R$ is a recursive relation
$\Pi_1$: $\{x : \forall y\, (\, R(x, y)\, )\}$      where $R$ is a recursive relation
$\Sigma_n$: $\{x : \exists y_n\, (\, P(x, y_1, ..., y_n)\, )\}$      where $P$ is a $\Pi_{n-1}$ relation
$\Pi_n$: $\{x : \forall y_n\, (\, S(x, y_1, ..., y_n)\, )\}$      where $S$ is a $\Sigma_{n-1}$ relation

From these, we define:
$\Delta_n$: $\Sigma_n \cap \Pi_n$

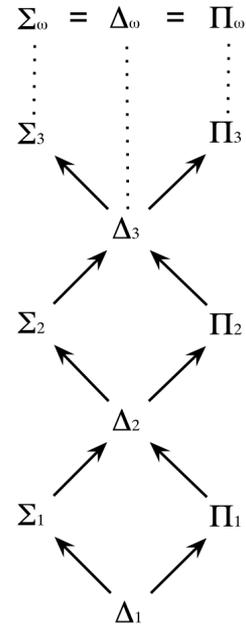

*Figure 8*

Figure 8 shows the large-scale structure of the arithmetical hierarchy, with the arrows representing strict set inclusion. $\Sigma_\omega$ represents the union of $\Sigma_n$ for all $n$ and similarly for $\Pi_\omega$ and $\Delta_\omega$.

Of particular note, $\Delta_1$ is the set of recursive relations, $\Sigma_1$ is the set of recursively enumerable relations and $\Pi_1$ is the set of co-recursively enumerable relations. We can see, for instance that the relation defined by the halting function is in $\Sigma_1$ by considering that the relation $HaltsBy(x, t)$ that is true if and only if the Turing machine and input encoded by $x$ halts within $t$ time-steps. $HaltsBy$ is clearly recursive and thus $\exists t\, (\, HaltsBy(x, t)\, )$ – which represents the halting relation – is in $\Sigma_1$. Similarly, we can show that the relation defined by the looping function can be represented as $\forall t\, (\, \neg HaltsBy(x, t)\, )$ and thus is in $\Pi_1$.

The usefulness of the arithmetical hierarchy is expanded by the following theorems:

*The representation theorem*:
     For any relation $R$, $R$ is in the arithmetical hierarchy iff $R$ is definable in elementary arithmetic

*The strong hierarchy theorem*:
     $R \in \Delta_{n+1}$ iff $R$ is recursive in $\varnothing^{(n)}$ (i.e. the $n^{th}$ jump of the empty set as discussed in section 3.1)
     $R \in \Sigma_{n+1}$ iff $R$ is recursively enumerable in $\varnothing^{(n)}$
     $R \in \Pi_{n+1}$ iff $R$ is co-recursively enumerable in $\varnothing^{(n)}$

## 5.2  Comparative Powers

The arithmetical hierarchy therefore gives us a generalisation of the concepts of recursiveness and recursive enumerability which allow us to measure and compare the power of the hypermachines that have been explored. This power is measured in terms of the predicates over **N** that can be computed by each hypermachines. This generalises to functions from **N** to $\{0,1\}$, sets over **N**, or computable reals.

$\Delta_1$:	Turing machines
	O-machines with $\varnothing$ oracle (or any recursive oracle)
	Turing machines with $\varnothing$ initial inscription (or any recursive inscription)
	asynchronous networks of Turing machines with recursive timing functions
	error prone Turing machines with recursive error functions
	coupled Turing machines with recursive coupled functions
	probabilistic Turing machines with recursive random choices

$\Sigma_1$:	fair non-deterministic Turing machines[19]
	accelerated Turing machines[20]

$\Delta_2$:	O-machines with $\varnothing'$ oracle
	Turing machines initially inscribed with $\varnothing'$
	infinite time Turing machines (in $\omega$ time-steps)[21]

$\Delta_{n+1}$:	O-machines with $\varnothing^{(n)}$ oracle
	Turing machines initially inscribed with $\varnothing^{(n)}$
	infinite time Turing machines (in $\omega \cdot n$ time-steps)

$\Delta_\omega$:	O-machines with $\varnothing^{(\omega)}$ oracle
	Turing machines initially inscribed with $\varnothing^{(\omega)}$
	infinite time Turing machines (in $<\omega^2$ time-steps)

There is also a large collection of machines that can compute all predicates over the natural numbers. This involves much more power than can be expressed in the arithmetical hierarchy (there are uncountably many predicates over **N**, but only a countable subset are expressible in the arithmetical hierarchy).

every predicate over **N**:
	O-machine (with arbitrary oracle)
	Turing machine (with arbitrary initial inscription)
	coupled Turing machines (with arbitrary coupled function)
	error prone Turing machines (with arbitrary error function)
	asynchronous networks of Turing machines (with arbitrary timing functions)
	probabilistic Turing machines (with arbitrary random events)

	infinite state Turing machine

---

[19] For a proof of this, see [45].

[20] For a proof of this, see the appendix.

[21] The results for the power of infinite time Turing machines with limited time-steps are lower bounds only. The upper bounds are not known, but I would be quite surprised if they are different. For the proofs, see the appendix.

Finally, while infinite time Turing machines can be shown to compute all predicates expressible in first order arithmetic in $<\omega^2$ time-steps, the set of predicates that they can decide in arbitrary time exceeds the arithmetical hierarchy. The exact power is unknown, but is clearly much less than *all* predicates over **N** (since there are only countably many infinite time Turing machines).

Some of these results have particularly interesting consequences. For instance, consider the models that can compute exactly the $\Sigma_1$ predicates. These models can compute the halting function, returning 1 if the program halts on its input and 0 otherwise, but cannot compute the looping function which would only require negating this result. This is quite an interesting result, because while accelerated Turing machines and fair non-deterministic Turing machines can compute both the halting function and the negation function ( $f(x) = 1 - x$ ), they cannot compute the composition of these functions. The functions computable by these models are thus not closed under composition, which is most unusual amongst serious models of computation. Since we would expect the functions computable by physical machines to closed under composition (the output of one machine can be used as input of another) and we know that negation is computable by physical machines, it seems highly implausible that $\Sigma_1$ is the limit of physical computation available in our universe.

The absence of closure under function composition also makes these models somewhat strange as theoretical devices. Indeed, it is probably best to imagine such devices as part of a larger machine. As Copeland [16] shows, we could embed one of these devices inside a Turing machine as an oracle. Since they can function as oracles for $\varnothing'$, the resulting O-machine could compute the entirety of $\Delta_2$.

It is also of special interest that some of these hypermachines can compute $\Delta_\omega$. This means that they can compute all the predicates that can be expressed as sentences in first order arithmetic. From the Representation Theorem, we can see that this includes all formulae of arithmetic that have free variables. What about those that do not have free variables? Consider an arbitrary formula of arithmetic that has no free variables. This sentence must include an = symbol because it is the only way of constructing an atomic formula in this language. We can construct a new formula by adding the term '+ $x$' (where $x$ is a new variable) to one side of the = sign. This gives us a new formula in arithmetic which expresses a predicate and is thus computable by this group of hypermachines. In particular, these hypermachines can compute the answer when $x$ is equal to 0. In such a case, this formula is true if and only if the original formula is true. Thus, hypermachines that can compute $\Delta_\omega$ can compute the truth of *all* sentences of arithmetic.

This is not a violation of Gödel's Incompleteness Theorem because his theorem only applies to formal systems (which are defined to have the same computational power as Turing machines). Gödel's methods cannot construct a true but unprovable sentence for these hypermachines because their workings cannot be captured in first order arithmetic. This result obviously has considerable implications for the common (mis)interpretation of Gödel's Incompleteness Theorem as stating that no conceivable machine can prove exactly the true sentences of arithmetic.

It is also important to consider the implications of those models that can compute all predicates over the natural numbers. Doing so requires an uncountable infinity of machines – at least one per predicate. Thus, the machines require infinite specifications (in most cases, this corresponds to specifying the non-classical part as well as the classical part). It is not at all clear how we should do so. While some infinite specifications are finitely definable in, say, first order arithmetic, each language with a finite alphabet has only a countable amount of finite strings. This is not enough to specify all the machines, so any attempt to finitely specify these machines through the use of a mathematical formalism must either use an infinite alphabet or infinitely many mathematical languages. Neither case looks any simpler than the original infinite specifications. It thus seems that while the universe may well have processes that behave according to these models, we will never be able to harness them to compute arbitrary functions.

Machines that can compute all predicates over **N** also seem to be of little theoretical interest. All of those examined here can compute a given predicate over **N** using a simple infinite lookup table (see the study of infinite state Turing machines in section 3.7 for a good example of this). They therefore do not need any algorithmic sophistication to compute these predicates. This is a similar situation to Turing machines computing predicates over a finite set – the machines are so powerful compared to the problem that only a degenerate kind of computation occurs.

# Chapter 6

# Implications

## 6.1 To Mathematics

One of the first uses of the Turing machine model was to show that the *Entscheidungsproblem* (the decision problem for the predicate calculus) is undecidable. By the Church-Turing Thesis, all functions computable by unaided mathematicians working without insight are computable by a Turing machine and since there is no Turing machine to solve the decision problem for predicate calculus, it was concluded that no unaided mathematician could solve it without insight.

The concept of hypercomputation sheds new light upon this and similar results. It prompts us to ask whether the Entscheidungsproblem could be solved by a mathematician that was aided by some kind of machine. The answer to this is not within the domain of mathematics. Instead it depends on the nature of physical reality. If it is possible to harness physical processes to build a hypermachine, then a mathematician may well be able to use this to decide whether arbitrary formulae of predicate calculus are tautologies.

Thus the claims that such problems are 'undecidable' or 'unsolvable' are misleading. As far as we know, in 100 years time these problems might be routinely solved using hypermachines. Mathematicians may type arbitrary diophantine equations into their computers and have them solved. Programmers may have the termination properties of their programs checked by some special software. We cannot rule out such possibilities with mathematical reasoning alone. Indeed, even the truth or otherwise of the Church-Turing Thesis has no bearing on these possibilities. The solvability of such problems is a matter for physics and not mathematics[22].

To show that a problem is physically undecidable should be a two part process. First, we mathematically determine how difficult the problem is using the techniques of recursion theory. Then we must ascertain whether or not it is physically possible to solve problems of that level of difficulty. Traditional arguments for undecidability have only focussed on the first of these steps, assuming that no non-recursive function is physically computable. This often unstated assumption is not only a very strong claim about the nature of our universe, but also lacking in physical evidence. It is not a mathematical claim and it is not appropriate to use it in a mathematical proof for anything more than speculation about the physical implications.

---

[22] A similar view is held by Deutsch, Ekert and Luppachini [24] who state that 'though the truths of logic and pure mathematics are objective and independent of any contingent facts or laws of nature, our *knowledge* of these truths depends entirely on our knowledge of the laws of physics'. See also Deutsch [23] for a further discussion of his view on the relationship of physics to mathematics.

This distinction that I am drawing between the mathematical and physical sides to an undecidability result generalises to Chaitin's results about randomness in mathematics.

Chaitin proved that $\Omega$ is 'random' and that its non-recursive information is so much more dense than $\tau$ that no formal system of size *n* bits can determine more than *n* bits of $\Omega$. He has also shown that analogous results occur in arithmetic, leading to systems of equations whose solutions do not have any recursive pattern. These results are nothing short of spectacular, but they are only results about the complexity of $\Omega$ relative to Turing machines and formal systems that are equally powerful. While it is sometimes claimed [4, 10] that $\Omega$ is completely patternless, this is not true. If it was, it would not be so easily describable. Indeed $\Omega$ is recursively enumerable and thus computable by every hypermachine presented in this paper. If a single one of these hypermachines can be built, we may well see people computing the digits of $\Omega$ as they now compute the digits of $\pi$.

Chaitin's results about $\Omega$ show not that some mathematical objects are patternless, but that they seem patternless relative to a certain amount of computational power[23]. The results do not show that some facts (such as the value of a particular bit of $\Omega$) are true for no reason, but that a certain level of computational power is insufficient to show why these facts are true. If hypercomputation is found to be physically impossible, then we could expand these results, with *our universe itself* being insufficient to find pattern in $\Omega$ or reasons for its digits being what they are. In any event though, it is not an absence of pattern, but rather a complexity of pattern, or an inaccessibility of pattern.

Gödel's results are similarly affected. While he successfully showed that no formal system can prove all truths of arithmetic, the implications that are commonly drawn from this are flawed. Even in the foreword of Davis' *The Undecidable* [20] we are told Gödel proved that 'systems of logic, no matter how powerful, could never admit proofs of all true assertions of arithmetic'. While this particular claim is perhaps more careless than misguided, it is indicative of the general presentation of Gödel's results. What he really showed was that no formal system *of equivalent power to the Turing machine* could enumerate all of the truths of arithmetic.

Many other systems of logic with non-recursive sets of axioms or rules of inference are possible. One could, for example, program an infinite time Turing machine or an O-machine with $\emptyset^{(\omega)}$ oracle to enumerate all the truths of arithmetic. It certainly requires a lot of computational power to prove all these truths with no falsehoods, but there are systems of logic which achieve this and it may even be physically possible to implement them. What we can mathematically prove depends on physics. In universes with different physics, different levels of computation would be possible and thus different sets of theorems would be provable from a given set of axioms[24].

---

[23] See section 6.3 for a more detailed development of this idea.

[24] See Deutsch [23] for a similar claim about the dependency of mathematical reasoning on physics.

And so each of these famous negative results does not go so far as is commonly argued. The limitations present in the results of Gödel, Turing and Chaitin are all relative to Turing machines. As far as we know, it may soon be routine to decide the truth of arbitrary statements of arithmetic, solve the halting problem and compute the digits of Ω. To find out, we need to complete the second part of the method mentioned above. The difficulty of these tasks has been determined – now we need to find out whether the universe has enough resources to compute them.

Even if it does, the ideas of Gödel, Turing and Chaitin are still very relevant. As well as showing the limits of Turing machines, they present a general method for finding problems that are unsolvable for a given type of machine. Presumably, we will find similar cases of incompleteness, undecidability and randomness relative to our most powerful implementable machines, whatever they happen to be. My argument is not with the limitation of our mathematical knowledge, but rather that we so far have no specific examples which we can say are incomplete, uncomputable or random relative to our physical laws.

## 6.2 To Physics

Computability is a very important physical issue. As argued above, the fundamental question regarding how much we can compute is a question of physics, and not mathematics as has been traditionally believed. We are now quite familiar with other physical limits such as the speed of light, the second law of thermodynamics or the conservation of energy. As we become more aware of the importance of computation at both fundamental and everyday levels, our lack of a physical principle explaining the limits of computation becomes more glaring.

Feynman [26] was the first to address such a gap in terms of the limits of computational efficiency, introducing the idea of quantum computing in order to close the apparent gap between what can be simulated efficiently and what can happen (efficiently) in quantum mechanics. His idea of harnessing these computationally complex results that nature was producing so quickly in quantum systems slowly built up the field of quantum computation which seriously addresses the issues of the universe's limits for the efficiency of computation.

We are perhaps seeing the beginning of a similar movement for exploring the physical limits of computational power. A few physicists have recently suggested possible ways that quantum mechanics and relativity may allow not only complexity gains, but computational gains. Hogarth [30, 31] has argued that general relativity may allow an infinite number of steps to be performed in a finite time. Calude and Pavlov [7] have shown how quantum mechanics may allow an infinite amount of steps to be performed in parallel, allowing a stochastic solution to non-recursive functions. Kieu [32] demonstrates a similar result.

There has also been some interest by physicists in exploring the fundamental connections between computation and physics. Deutsch [22, 23, 24] explores such connections, speculating on the limits of computation and addressing the fundamental question of how difficult it is to simulate the universe. He proposes the 'Church-Turing Principle': that the universe is capable of containing a universal machine that can simulate arbitrary systems of the universe.

In chapter 2, I presented my own physical theses about the physical limitations of the universe, suggesting ways in which the universe may not be self-simulable due to an absence of harnessable hypercomputation. I also made a distinction between mathematical and physical harnessability, showing how it may be possible for the universe to be unsimulable through a mathematical representation, but still simulable. A universe that is classically simulable, but for the existence of a number of rods of the same non-recursive length is such an example. We may be able to use one of these rods to simulate the others even if we cannot produce a finite mathematical description of its length. Such possibilities present new challenges for the foundations of physics and the philosophy of science.

## 6.3 To Computer Science

Hypercomputation clearly has considerable implications for computer science, yet it has received almost no attention there. Part of the reason for this is that few people realise that it is logically possible to exceed the power of the Turing machine. It is often believed that the Church-Turing Thesis rules out the possibility of hypercomputation. Another common misconception is that the proof of the uncomputability of the halting function by Turing machines makes it paradoxical for *any* machine to compute it. Finally, people dismiss the possibility of hypercomputation as completely unrealistic or in violation of some physical laws.

The first two of these reasons for ignoring the possibility of hypercomputation are simply false. The third would have some force if it was used with an understanding of the many different resources that could be used to achieve hypercomputation, but there has been no detailed analysis of whether these resources are physically realisable. It is thus very unfortunate that the reasons above have caused so many people to ignore the fascinating possibilities that hypercomputation offers.

One of the main reasons that people still hold these misconceptions is that the basic issues of hypercomputation – particularly those three just mentioned – are often not well understood and consequently the possibility of hypercomputation is not presented in textbooks or lectures. While there are exceptions to this (most notably brief presentations of O-machines), students still tend to finish their computer science education with a firm belief that the Turing machine has limits on its computation, but that it is as powerful (in terms of computability) as any other possible machine.

There are very important issues regarding computability at stake here. What forms can an algorithm take? What is computable in our universe? These are questions that one would expect to attract many researchers. They are unanswered (and relatively unexplored) questions that are of great importance to computer science.

One example of the gains that can be made by applying the lessons of hypercomputation to classical fields is in algorithmic information theory. Hypercomputation provides a very natural relativisation of algorithmic information theory. Previously, we looked at algorithmic information theory in the context of a person, Alice, sending compressed messages to a friend, Betty, who uses a Turing machine to decode them. Of particular interest was the case of sending infinite messages. For some messages (such as the digits of π) there is a finite program that can be sent instead. For others (such as the digits of τ) there is no finite program that can be sent which generates this sequence. Thus if Alice somehow has access to τ, she needs to send an infinite number of bits to Betty to share this information with her. It can thus be said that τ contains an infinite amount of algorithmic information.

As was discussed in chapter 1, for Betty to receive $n$ bits of τ, Alice only needs to send $log(n)$ bits of information because τ's information is spread very thinly. For Ω on the other hand, Alice needs to send $n$ bits for Betty to determine $n$ bits, thus Ω is said to be random.

In a *relativised algorithmic information theory*, we can provide Betty with different types of computer. For example, consider a faulty Turing machine that can only move its tape head to the right. It can be shown that such a machine can compute all the rationals, but no irrational numbers. Therefore, if Alice and Betty were restricted to such machines and Alice wanted to send the number π to Betty, she would be forced to send an infinite bitstring. With this model of computation, π has an infinite amount of algorithmic information.

We can also look at what happens when Betty has access to a hypermachine. If Betty has access to, say, an O-machine with the halting function oracle, then Alice can send her a finite program for computing τ or Ω. Relative to an O-machine with the halting function oracle, τ and Ω both have finite algorithmic information. However, the bits that describe this machine's halting function (or the truths of arithmetic) would still require an infinite message.

With a relativised algorithmic information theory, we would say that τ contains an infinite amount of information relative to a Turing machine, that Ω is random relative to a Turing machine and that π contains an infinite amount of information relative to a 'right-moving Turing machine'.

This appears to be a deep and natural generalisation from the standard theory and is not limited by any reliance on the actual computational power of our universe. Instead, it provides a concept of information that is relative to the level of power available and thus relative to our particular laws of physics. By exploring the ideas involved in hypercomputation, we should be able to generalise and extend the classical results of other areas of computer science in a similar manner. A study of complexity theory over these machines, for instance, would no doubt be very interesting and illuminating.

Finally, if we can actually build any of these hypermachines, there will obviously be large implications for computer science. Hypercomputation would presumably become mainstream computer science as people move to take advantage of these new powers. It is worth noting that the uptake of hypercomputation could still be limited by such things as impractical construction, price, or even speed of computation – it may be that the only buildable hypermachines are exponentially slower than Turing machines. In any event, there is no doubt that if we discover that hypercomputation is physically possible and harnessable, its implications to computer science will be immense.

## 6.4 To Philosophy

The possibility of hypercomputation also bears important implications in several branches of philosophy. One of these is the philosophy of mind. The computational theory of the mind looks at the possibility that the computations performed in the brain give rise to the mind. It is often claimed that if the mind is indeed computational in nature, then we could simulate the brain (and mind) on a Turing machine. This argument often appeals to the Church-Turing Thesis to show that all computation can be simulated on a Turing machine. This mistaken belief in the power of the Church-Turing Thesis pervades the literature and is discussed in great depth by Copeland [13]. Copeland also shows how the claims made in the computational theory of mind can be generalised to allow for the possibility of hypercomputational processes in our brains giving rise to our minds [15].

Hypercomputation also has significant implications in epistemology. In classical computation, there is a fixed limit to the conclusions that can be validly generated from a set of premises, but hypercomputation generalises this with the amount of valid inferences growing with the level of computational resources available. This not only relativises provability and inference, but makes the actual amount of inference that is possible in the real world dependent upon the true nature of physics and not *a priori* as it is generally believed.

In addition, epistemology involves assessing the complexity of different ideas. Simplicity is a crucial notion when deciding which of two theories to prefer and thus the changes to algorithmic information theory brought about by hypercomputation become relevant. To measure the complexity of a theory, we must do so relative to a given level of computational power and thus an understanding of hypercomputation seems necessary here.

It is not clear how large the impact of hypercomputation in these areas will be, but it should not be overlooked. Several areas of philosophy involve computational concepts and the widespread lack of knowledge about hypercomputation will no doubt have caused important possibilities to be overlooked. Hypercomputation is very relevant to philosophy and must be accepted as a logical possibility that needs careful consideration and investigation.

## Chapter 7

# Conclusions and Further Work

### 7.1 Conclusions

Models of computation that compute more than the Turing machine are conceptually unproblematic and may even be buildable. This has many important implications. The most direct of these is that mathematical proofs that certain functions are 'uncomputable' must be more clearly understood as shorthand for 'uncomputable by a Turing machine' and not as a limitation on mathematics. As I have argued in chapters 5 and 6, the proofs of incompleteness and randomness in arithmetic must also be seen as merely relative to Turing machines. By opening the possibility of physical machines that compute the halting function, Ω and all the truths of arithmetic, hypercomputation shows us that the great negative results of 20[th] Century mathematics are not so negative as it was first believed.

In doing so, hypercomputation follows its analogue in complexity theory – quantum computation – in elucidating the connections between mathematics, physics and computer science. As quantum computation showed that mathematics alone cannot provide the limits of physically practical computation, so hypercomputation shows that mathematics alone cannot provide the outer limit for physically possible computation.

In chapter 2, I looked at three different ways in which physical computation could be limited. It may well be that there are physical processes that, while not simulable by a Turing machine, are not harnessable for the computation of any particular non-recursive functions. Furthermore, it may be that there exist non-recursive physical processes that cannot be simulated with the resources available within the universe. The logical independence of the computational limits on physical processes, physically harnessable processes and mathematically harnessable processes allows many different ways in which our universe could be hypercomputational – each with its own implications for physics and the philosophy of science.

With the increasing understanding of the misconceptions surrounding the Church-Turing Thesis, prematurely closed possibilities – such as a hypercomputational theory of mind – can be explored. Through non-classical models of computation, such as infinite time Turing machines, we can also attempt to gain a deeper understanding into the nature of algorithms. Regardless of the actual computational limits of our universe, I have no doubt that the study of hypercomputation will lead to many important theoretical results across computer science, philosophy, mathematics and physics.

## 7.2 Further Work

The ideas explored in this report are only the very beginnings of a comprehensive theory of hypercomputation and there are thus many avenues for further work.

Theses **B**, **C** and **D**, in chapter 2 provide interesting insights into different computational limits in physics. However, they rely on the terms 'simulable', 'mathematically harnessable' and 'physically harnessable'. While I have tried to give useful definitions of these terms, more work is required to create precise concepts with which to assess the availability and usefulness of hypercomputation in our universe. This is quite related to the problem of defining what it means for a physical process to compute a given function and – while there are some difficulties with this – the techniques used by Copeland [12] and Deutsch [22, 23] look particularly promising.

While the theses I have provided present an interesting stratification of the computational powers of different possible universes, they could be extended to provide a more fine-grained classification. In particular, there would be considerable benefit in extending them to assess the universe's simulability by arbitrary models of computation instead of just Turing machines. It would also be beneficial to explore the relationship between my theses and Deutsch's Church-Turing Principle.

In this report, I have not had the space to present every hypermachine that has been explored. In particular, I have not been able to look at non-uniform circuit families, type-2 machines [56], Hava Siegelmann's processor networks [43] or Lenore Blum, Mike Shub and Steve Smale's machines for computation on arbitrary rings (BSS machines) [2]. The presentation of these models and the linking of them into the framework of this paper would obviously be of great use.

In chapter 4, I looked briefly at a way of exploring models between the power of an accelerated Turing machine and an infinite time Turing machine. I have given only a few proofs and results about such machines in this paper and a formal presentation of such machines and their capabilities would be most interesting.

It would also be of considerable interest to examine the more 'realistic' resources of interaction and infinite run-time as discussed by Jan van Leeuwen and Jirí Wiedermann [55]. They present good arguments for requiring a new model for the study of such machines due to their increased power over Turing machines. Looking at this in the context of hypercomputation may provide some useful results for the study of these systems.

In studying the capabilities of the hypermachines in this report, my use of the arithmetical hierarchy restricted my analysis to functions from **N** to {0,1} (or equivalently, reals or sets of naturals). It would be very nice to be able to extend this analysis to functions from **N** to **N**, providing a complete account of the capabilities of these machines. It would also be useful to discuss the capabilities of various machines that compute over **R**, including infinite time Turing machines, type-2 machines and BSS machines.

I have also only been able to give a brief introduction to the relativisation of algorithmic information theory to arbitrary computers. This is a very promising area which is practically unexplored. An in-depth look at the generalised notions of incompressibility and randomness would provide considerable interest.

While I have shown that there are generalised algorithms for enumerating the truths of arithmetic, computing the halting function and computing $\Omega$, there is a strong sense in which we remain under the shadow of incompleteness, uncomputability and randomness. It seems as though the methods used by Gödel, Turing and Chaitin to show these limitations of mathematics can be extended to models beyond the Turing machine. Indeed it seems quite likely that while these famous arguments do not establish any absolute limits for mathematics, they can be used to construct similar limitations for any particular model of computation. Whether or not this can be achieved, and the exact form that such generalised arguments would take is of considerable interest and is certainly a key area for further research.

# Appendix

Let 𝒫(M) be the set of predicates over **N**, decided by the model M.

### 𝒫(accelerated Turing machine) = $\Sigma_1$

Every predicate in $\Sigma_1$ can be expressed as $\exists y\, ( R(x, y) )$, where $R$ is a recursive predicate. An accelerated Turing machine, $A$, given $x$ as input can compute $R(x, 1)$ in a finite number of time-steps, marking its designated square if $R(x, 1)$ is true. If $R(x, 1)$ is not true, $A$ can proceed to compute $R(x, 2)$, $R(x, 3)$ and so on, marking the designated square and halting if it finds an instance which is true. After two external time units, the designated square will contain **1** if and only if $\exists y\, ( R(x, y) )$. Thus, every predicate in $\Sigma_1$ has an equivalent accelerated Turing machine.

Therefore, 𝒫(accelerated Turing machine) $\supseteq \Sigma_1$

Given an accelerated Turing machine $A$, let $R_A(x, y)$ be true in and only if when $A$ acts on $x$, it marks the designated square within $y$ time-steps. $\exists y\, ( R_A(x, y) )$ is true if and only if $A$ returns 1 when acting on $x$. $\exists y\, ( R_A(x, y) )$ is false if and only if $A$ returns 0 when acting on $x$. Since $R_A(x, y)$ is recursive for all $A$, $\exists y\, ( R_A(x, y) )$ is in $\Sigma_1$ for all $A$. Thus every accelerated Turing machine has an equivalent predicate in $\Sigma_1$.

Therefore, 𝒫(accelerated Turing machine) $\subseteq \Sigma_1$

Therefore, 𝒫(accelerated Turing machine) $= \Sigma_1$

### 𝔓(ω-time Turing machine) ⊇ Δ₂

By the Strong Hierarchy Theorem, a predicate is Δ₂ iff it is recursive in ∅'. Thus, for ω-time Turing machines to compute Δ₂, it suffices to simulate each halting O-machine with oracle ∅', that computes a predicate over **N** (i.e. a function from **N** to {0,1}). Let $O$ be an arbitrary halting O-machine with oracle ∅' and let $W$ be an ω-time Turing machine. $W$ can perform a non-deterministic simulation of $O$, branching into two computations whenever $O$ makes a request to its oracle. One computation branch assumes that $n$ was in the oracle set and one assumes it was not. Let us label these branches 1 and 0. Each branch is to be simulated in parallel (see figure 9 for a representation of $W$'s computation tree).

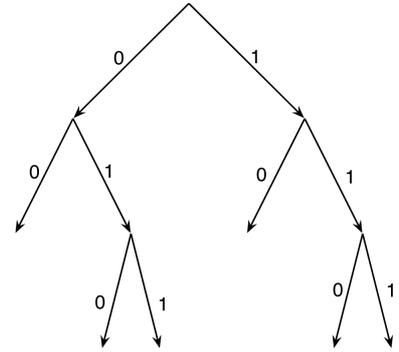

*Figure 9*

Since ∅' is the set of all encodings of Turing machines and input that halt, when $O$ asks its oracle if $n$ is in the oracle's set, this is equivalent to asking if $n$ represents a halting Turing machine/input combination. Every time a $O$ makes a request to its oracle, $W$ also begins a parallel simulation of the appropriate Turing machine on its input. If the Turing machine halts on its input, it is clear that $n$ was in the oracle set, so $W$ ceases its simulation of the computation sub-tree where $n$ was assumed to not be in the set. Finally, every time the *leftmost* computation branch (that has not been ceased according to the previous rule) halts, its returned value (either a 0 or a 1) is placed on the first square of $W$'s tape (a place specifically reserved for this).

Since $O$ halts, it can only make a finite number of requests to its oracle (the other branches that are simulated which don't reflect the advice that $O$ would actually receive from its oracle may not halt, but this need not concern us). Therefore, after a finite number of time-steps, all the oracle requests actually made by $O$ that return 'true' (or 1) have been identified as such and the respective computation sub-trees in which an answer of 'false' (or 0) was being simulated have been ceased and removed – of course $W$ can not determined when this happens, but it is only important that it will eventually happen. Also, if there are still two computation branches for an oracle request, the left branch (the 0 branch) corresponds to the oracle's answer.

Therefore, the computation that corresponds to the leftmost branch in this computation tree corresponds to the computation of $O$ (and halts after a finite time). Since this branch has halted and all the 'guessed' oracle answers in this branch have stabilised, this branch will remain the leftmost branch. Therefore, the result of $O$'s computation will be stored in the first square of $W$'s tape from this point on. Thus, $W$ computes the same function as $O$.

Therefore, 𝔓(ω-time Turing machine) ⊇ Δ₂

### 𝒫(ω·*n*-time Turing machine) ⊇ Δ*ₙ*₊₁

In ω time-steps, an infinite time Turing machine can simulate every Turing machine / input combination, putting a 1 on the appropriate square if and when the Turing machine / input combination represented by *n* halts. After ω time-steps, the tape will contain $\varnothing'$ (and can still have an infinite amount of room for scratch working). In another ω time-steps, the machine can do the same, but simulating every O-machine with $\varnothing'$ as its oracle, producing $\varnothing''$ on its tape. This can be repeated, producing $\varnothing^{(n-1)}$ in ω·(*n*-1) time-steps.

By the Strong Hierarchy Theorem, a predicate is $\Delta_{n+1}$ iff it is recursive in $\varnothing^{(n)}$. Therefore, once the machine has $\varnothing^{(n-1)}$ on its tape, it has ω time-steps remaining in which to simulate an O-machine with $\varnothing^{(n)}$ as its oracle. This can be accomplished by trivially adapting the argument given on the previous page for using ω time-steps to simulate an O-machine with $\varnothing'$ as its oracle.

Therefore, 𝒫(ω·*n*-time Turing machine) ⊇ $\Delta_{n+1}$

### 𝒫(<ω²-time Turing machine) ⊇ $\Delta_\omega$

Any predicate in $\Delta_\omega$ is in $\Delta_n$ for some *n*. Therefore any predicate in $\Delta_\omega$ is decided by some infinite time Turing machine in ω·*n* time-steps. Thus each predicate in $\Delta_\omega$ is decided by some infinite time Turing machine in less than ω² time-steps.

Therefore, 𝒫(<ω²-time Turing machine) ⊇ $\Delta_\omega$